\documentclass[aop,citesort,MSNbibl,noautosecdot,dvips]{arximspdf}
\usepackage{graphicx}

\doi{10.1214/10-AOP618}
\volume{39}
\issue{5}
\pubyear{2011}
\firstpage{1864}
\lastpage{1895}

\makeatletter

\newcommand{\nicefrac}[2]{{#1}/{#2}}
\newcommand{\dotsb}{\cdots}

\newtheorem{theorem}{Theorem}
\newtheorem{lem}{Lemma}
\newproclaim{defn}{Definition}

\newcommand{\dotcup}{\mathop{\ooalign{\hfil$\bigcup$\hfil\cr\hfil$\cdot$\hfil\cr}}}
\newcommand{\chem}{\operatorname{int}}

\makeatother

\begin{document}
\begin{frontmatter}

\title{Percolation on a product of two trees\thanksref{T1}}
\runtitle{Percolation on a product of two trees}

\thankstext{T1}{Supported in part by the Israel Science Foundation.}

\begin{aug}
\author[A]{\fnms{Gady} \snm{Kozma}\corref{}\ead[label=e1]{gadyk@weizmann.ac.il}}
\runauthor{G. Kozma}
\affiliation{Weizmann Institute of Science}
\address[A]{Department of Mathematics\\
Weizmann Institute of Science\\
Rehovot 76100\\
Israel\\
\printead{e1}} 
\end{aug}

\received{\smonth{3} \syear{2010}}
\revised{\smonth{10} \syear{2010}}

%
\begin{abstract}
We show that critical percolation on a product of two regular trees
of \mbox{degree $\ge$ 3} satisfies the triangle condition. The proof does
not examine the degrees of vertices and is not ``perturbative''
in any sense. It relies on an unpublished lemma of Oded Schramm.
\end{abstract}

%

\end{frontmatter}

\section{\texorpdfstring{Introduction.}{Introduction}}\label{sec1}

\subsection{\texorpdfstring{Schramm's lemma.}{Schramm's lemma}}\label{sub:Oded}

In 1998, while working on the exponential decay of correlations problem
(a problem which is still open), Oded Schramm proved a lemma which
solved the problem in an averaged sense. The lemma generated a lot
of excitement in the community of researchers of percolation on groups
at the time, and there was hope that it would lead to a full solution
of the exponential decay problem, the $p_{c}<p_{u}$ problem and other
related problems on nonamenable Cayley graphs. That hope never materialized,
Schramm never published the lemma (you will see it mentioned here
and there in papers of the period, e.g., in the last paragraph of
\cite{BLPS99}) and moved on to other topics.

Let us describe the settings of Schramm's lemma in its original formulation.
We are given a nonamenable Cayley graph%
\setcounter{footnote}{1}\footnote{Cayley graphs, amenable groups, and any other notion not explicitly
defined may be found in Wikipedia.%
} $G$. Denote the identity element of the group by $0$ (by which
I definitely do not insinuate that the group is Abelian) and for any
$x\in G$ let $|x|$ denote the graph distance between $0$ and $x$
(which is the same as the distance in the word metric on the group
with the given generators). Let $p_{c}$ be the critical probability
for percolation on $G$. The exponential decay of correlations problem
is the conjecture that at $p_{c}$ one has that the connection probability
$\mathbb{P}(0\leftrightarrow x)$ decays exponentially in $|x|$.
Schramm's lemma states that this is true in the following interesting
averaged sense. Take random walk $R$ on the original Cayley graph
$G$. Then
%
\begin{equation}\label{eq:oded orig}
\mathbb{P}\bigl(0\leftrightarrow R(n)\bigr)\le c^{n}.
\end{equation}
Let us stress again that the random walk is on the original graph
$G$ and is only some way to average $\mathbb{P}(0\leftrightarrow x)$.
The probability in (\ref{eq:oded orig}) is over both the walk and
the percolation.

We will now sketch Schramm's argument. If the sketch is too dense
do not despair---a very close result, Lemma \ref{lem:Oded}, will
be presented below, page~\pageref{lem:Oded}, with all details (perhaps
too many details). For general background on percolation, see the book
\cite{G99}. For percolation, random walk and branching processes
on transitive graphs see the book \cite{LP}, especially Chapters
7 and~8. Finally, note that all notations and conventions used in this
paper are collected in Section \ref{sub:Notations}, page \pageref
{sub:Notations},
for convenience.

\subsubsection*{\texorpdfstring{Proof sketch.}{Proof sketch}}
Examine the following process on our graph $G$. Fix
some~$m$. We start with $m+1$ ``particles'' at $0$. Each particle
does $n$ steps of simple random walk on $G$, independently of all
other particles. After $n$ steps, it divides into $m$ particles,
and each one does $n$ additional steps, divides and so on. In other
words, a more-or-less standard branching random walk. The particles
can be mapped to the edges of an~$m+1$ regular tree $\mathbf{T}$.
The original~$m+1$ particles are mapped to the edges coming out of
the root of $\mathbf{T}$. Then for each of these particles, say it
corresponds to the edge $(\mathbf{r},\mathbf{v})$ with $\mathbf{r}$
being the root of the tree and $\mathbf{v}$ one of its children,
map the $m$ descendents of the particle to the $m$ edges coming
out of $\mathbf{v}$ not in the direction of the root. Continue inductively.

We now throw the percolation into the mix as follows. We define a
random subgraph $\mathbf{W}$ of our tree $\mathbf{T}$ as follows.
Examine an edge $(\mathbf{v},\mathbf{w})$ of the tree. It corresponds
to some particle in the branching process which started from some
$v$, walked $n$ steps and landed on $w$. We declare that $(\mathbf
{v},\mathbf{w})\in\mathbf{W}$
if and only if $v\leftrightarrow w$. At the formal level, $\mathbf{W}$
is a random subgraph of $\mathbf{T}$ which is a deterministic function
of both the branching process and the percolation.

Note that $\mathbf{W}$ is, considered as a measure over subgraphs
of $\mathbf{T}$, invariant to the automorphisms of $\mathbf{T}$.
Please reflect on this fact for a minute as it utilizes a number of
features of the construction, with a particular emphasis on automorphisms
of $\mathbf{T}$ which do not preserve the root $\mathbf{r}$. First,
notice that it is crucial to start with $m+1$ particles at time $0$,
but split to only $m$ particles at each subsequent time. Next, note
that we used the reversibility of random walk: let $\varphi$ be an
automorphism which sends, say, some $\mathbf{v}$ which is a child
of $\mathbf{r}$ into $\mathbf{r}$. For $\mathbf{W}$ to be invariant
to $\varphi$ it is necessary that getting from the child to the father
has the same distribution as getting from the father to the child.
So our process must be time-reversible.

It is now time to fix $m$. We fix it such that the branching process
is transient, that is, such that with probability $1$ only a finite
number of particles ever return to~$0$. With this definition of $m$,
we get that the configuration $\mathbf{W}$ contains only finite components.
Indeed, by \cite{BLPS99} the percolation clusters are all finite.
Any cluster $\mathbf{C}$ of $\mathbf{W}$ is (a subset of) all returns
of the branching process to some percolation cluster, that is, a finite
collection of finite sets, so $\mathbf{C}$ is finite.

At this point, we apply the mass transport principle to the tree
$\mathbf{T}$
(which is of course itself the Cayley graph of a nonamenable group)
and we get that the expected degree of $\mathbf{r}$ is equal to the
expected average degree of the cluster of $\mathbf{r}$. But the cluster
of $\mathbf{r}$ is a finite subgraph of $\mathbf{T}$ so it is a
finite tree and for every finite tree the average degree is $\mbox{$<$}2$!
Hence, we get that the expected degree of $\mathbf{r}$ in the configuration
$\mathbf{W}$ is $\mbox{$<$}2$.

But what is the expected degree of $\mathbf{r}$ in $\mathbf{W}$?
It is exactly $m\cdot\mathbb{P}(0\leftrightarrow R(n))$! We get
\[
\mathbb{P}\bigl(0\leftrightarrow R(n)\bigr)<\frac{2}{m}
\]
so Schramm's argument terminates with the observation that if $m<\lambda^{-n}$
with $\lambda$ being the spectral radius of the random walk then
the branching process is transient. We get that $\mathbb
{P}(0\leftrightarrow R(n))<2\lambda^{n}$
which decays exponentially as for any nonamenable Cayley graph $\lambda<1$.

\subsubsection*{\texorpdfstring{Notes.}{Notes}}
\begin{longlist}[(1)]
\item[(1)] In fact, the process is transient even when $m=\lambda^{-n}$ for
$\lambda$ being the spectral radius. This property is called ``$\lambda
$-transience''
and holds for all nonamenable Cayley graphs. See \cite{W00}, Section II/7B.
\item[(2)] FKG implies that the function $f(n)=\mathbb{P}(0\leftrightarrow R(n))$
is supermultiplicative, that is, $f(a+b)\ge f(a)f(b)$. Therefore, if
we know that $f(n)<2\lambda^{n}$ then we actually know that $f(n)\le
\lambda^{n}$,
that is, the 2 may be dropped.
\item[(3)] A topic of some interest would be to generalize this argument from
a Cayley graph to a general transitive (nonamenable) graph. One
obvious point in the proof which would not apply to a general transitive
graph is the application of \cite{BLPS99}. We remark
that \cite{BLPS99} holds in larger generality than Cayley graphs,
it holds for any unimodular transitive graph. See \cite{BLPS99} or
\cite{LP}, Section~8.2, for details on unimodularity. Unimodularity or a
similar requirement might also be necessary if one wants to generalize
the claim that $\mathbf{W}$ is invariant to the automorphisms of
$\mathbf{T}$.
\end{longlist}

\subsection{\texorpdfstring{The triangle condition.}{The triangle condition}}

We say that a transitive graph (which may be amenable) satisfies the
triangle condition at some $p$ if
\[
\nabla_{p}:=\sum_{x,y}\mathbb{P}_{p}(0\leftrightarrow x)\mathbb
{P}_{p}(x\leftrightarrow y)\mathbb{P}_{p}(0\leftrightarrow y)<\infty.
\]
We will usually be interested in behavior at the critical $p$, and
denote $\nabla=\nabla_{p_{c}}$. The triangle condition was suggested
in 1984 by Aizenman and Newman \cite{AN84} as a marker for ``mean-field
behavior,'' a term from statistical physics which in our context
means that various quantities behave at or near $p_{c}$ as they would
on a regular tree. In particular, Aizenman and Newman proved that under
the triangle condition $\nabla<\infty$ one has that $\mathbb
{E}_{p}|\mathcal{C}(0)|\approx(p_{c}-p)^{-1}$
as $p$ tends to $p_{c}$ from below. Here, $\mathcal{C}(0)=\{
x\dvtx0\leftrightarrow x\}$
is the cluster of $0$, and $|\mathcal{C}(0)|$ is its size. See \cite{BA91}
for a proof that (again under $\nabla<\infty$) $\mathbb
{P}_{p_{c}}(|\mathcal{C}(0)|>n)\approx n^{-1/2}$
and $\mathbb{P}_{p}(|\mathcal{C}(0)|=\infty)\approx(p-p_{c})^{+}$.
See \cite{N87} for gap exponents and \cite{KN09} for intrinsic (a.k.a.
``chemical'')
exponents and the behavior of random walk on large critical clusters.
In short, under the triangle condition we have a very fine picture
of the behavior at and near criticality.

It is conjectured that the triangle condition holds in great generality.
A~``folk'' conjecture suggests that it holds for every transitive
graph for which the random walk triangle condition $\sum
G(0,x)G(x,y)G(0,y)<\infty$
holds, whe\-re~$G$ is the random walk Green function. I was not able
to find a reference for this precise formulation but a weaker one
is in \cite{S01}, Conjecture 1.2. Progress on this has been slow.
The most spectacular result is that the triangle condition holds for
the Cayley graph of $\mathbb{Z}^{d}$ with $d>6$ if one takes sufficiently
many generators; or with the standard generators $\pm e_{i}$ if $d$
is sufficiently large. This was achieved by Hara and Slade \cite{HS90}
using a technology known as ``lace expansion.'' See \cite{HHS08}
for generalizations to long-range models. The lace expansion is a
perturbative technique and generally requires the number of generators
be large. A different artifact of the perturbative nature is that
it seems that lace expansion is not suitable to show that $\nabla<\infty$
unless in fact $\nabla$ is quite close to $1$ ($\nabla\ge1$ always
due to the term $x=y=0$).

Going beyond $\mathbb{Z}^{d}$ there are two papers I am aware of
that establish the triangle condition. The first, \cite{W93}, establishes
it for $T\times\mathbb{Z}$ where $T$ is a~regular tree with degree
$\ge5$. The proof, very roughly, utilizes the fact that for a~tree
of degree $d$, $p_{c}=1/(d-1)$ but $\nabla_{p}<\infty$ for all
$p<1/\sqrt{d-1}$. This allows to make a relatively rough estimate
of $\nabla$ by path counting. A far more general result was achieved
in \cite{S01} which showed (among other things) that for any nonamenable
group, if one takes sufficiently many generators then the resulting
Cayley graph satisfies the triangle condition. Both results are significantly
easier than using lace expansion.

We may now state our result.
\begin{theorem}\label{theo1}
Let $T$ be a regular tree of degree $\ge3$. Then the product graph
$T\times T$ satisfies the triangle condition at $p_{c}$.
\end{theorem}

One may of course wonder whether one can achieve this result from
the ``highly nonamenable'' condition of \cite{S01}, but an inspection
shows that this would require $T$ to have degree $\ge7$. Also, the
approach taken here can be used also for ``stretched'' trees,
namely, suppose one replaces every edge with a path of length 100.
Naturally the resulting graph is no longer transitive (the degree
of vertices in the product of two stretched trees may be $2d$, $d+2$
or 4), but it is quasi-transitive, that is, its group of automorphisms
acts with finite orbits, and the argument presented below will work,
mutatis mutandis. Such a graph can have arbitrarily low Cheeger constant.
Similarly, it may have $\nabla$ arbitrarily large. It is probably
also worth comparing to~\cite{S02}, which shows for \textit{planar}
nonamenable transitive graphs, not quite the triangle condition,
but many mean-field exponents. This result is also nonperturbative
but relies critically on the planarity. It is easy to check that
$T\times T$
is not planar except when both have degree $2$.

The important property of $T\times T$ which is used is the large
number of symmetries. Examine the sphere $\{x\in T\times T\dvtx|x|=r\}$.
Clearly, $\mathbb{P}(0\leftrightarrow x)$ depends only on $|x_{1}|$
and $|x_{2}|$ where $x_{i}$ are the projections of $x$ to the two
trees. Thus, the sphere breaks down into $r+1$ of classes each of
which has exponentially many vertices. This is the crucial property.
Thus, the argument works for a product of two trees with different
degrees, or for products of three trees or more. It is probably possible
to formulate the result abstractly in terms of symmetries of the graph,
but I did not have a good formulation or a second interesting example,
so we will restrict ourselves to the simplest nontrivial example:
$T\times T$. In addition, let us remark that I could not make the
argument work for $T\times\mathbb{Z}$. Even though most vertices
do have many ``clones,'' some (namely, those on the same copy
of $\mathbb{Z}$ as $0$) have only one clone, and this is enough
to break the argument in its current form. We will return to this
topic in Section \ref{sub:TtimesZ}.

\subsection{\texorpdfstring{What has Schramm's lemma to do with the triangle
condition\textup{?}}{What has Schramm's lemma to do with the triangle
condition}}\label{sub:What}

Let us go back to the proof of Schramm's lemma. As already remarked,
the~ran\-dom walk is just some way to average the connection probability.
One may take any averaging method as long as it is time-reversible.
Taking a~3-regu\-lar tree $T$ as an example, one may replace the random
walk of length $n$~in Schramm's original argument with simply choosing
a random element of~dis\-tance $n$ from the root. One gets that
\[
\mathbb{P}(0\leftrightarrow R)<\frac{2}{m},
\]
where $R$ is a random element of distance $n$ from the root, and
$m$ needs to satisfy that the corresponding branching process is
transient. Of course, by the symmetries of the tree all elements of
distance $n$ from the root have the same connection probability so
the ``averaging'' performed by taking a~ran\-dom element has no
effect.

We will calculate the largest $m$ one may take below (it is the
claim in the proof of Lemma \ref{lem:Oded}), but we note for now
that the number of particles which return to $0$ after two steps
of the branching process is $m^{2}/2^{n}$. Therefore, it is reasonable
to assume that $m\simeq\sqrt{2^{n}}$ is the threshold for the branching
process to be transient, so we should have $\mathbb{P}(0\leftrightarrow
x)\lesssim2^{-|x|/2}$
and indeed a~slightly more precise calculation shows (still for the
tree) that
\[
\mathbb{P}(0\leftrightarrow x)<C|x|2^{-|x|/2}.
\]
This can now be summed explicitly to give
%
\begin{equation}\label{eq:triang in ball}
\sum_{|x|,|y|\le r}\mathbb{P}(0\leftrightarrow x)\mathbb
{P}(x\leftrightarrow y)\mathbb{P}(0\leftrightarrow y)\le Cr^{8}.
\end{equation}
The $8$ is of course not important---what is important is that
the sum grows only polynomially, even though the ball of radius $r$
grows exponentially. In other words, Schramm's lemma ``almost''
gives the triangle condition. This is the crucial observation on which
we rely.

\subsection{\texorpdfstring{Proof sketch.}{Proof sketch}}

To show the triangle condition, it is enough to show that $\mathbb
{P}(0\leftrightarrow x)\le C2^{-|x|(1/2+\varepsilon)}$
as then one can sum these explicitly. By the symmetries of our graph,
if we show that
\[
\mathbb{E}(\{x\dvtx0\leftrightarrow x,|x|\le r\})\le C2^{r(1/2-\varepsilon)}
\]
we will be done. The set $\{x\dvtx0\leftrightarrow x,|x|\le r\}$ is known
as the ``ball in the extrinsic metric''---meaning that we measure
the distance to $x$ by the distance inherited from the surrounding
graph, $|x|$.

If there is anything metamathematical to be learned from comparing
Section~3.2 in \cite{KN09} to \cite{KN}, it is that it is easier to work
with the intrinsic distance. In other words, rather than looking at
$|x|$, we should look at $d_{\chem}(x,y)$, the length of the shortest
path of open edges between $x$ and $y$. This idea has a checkered
past---in two dimensions the behavior of $d_{\chem}(x,y)$ is still
wide open---but in mean-field setting it has proved to be a useful
tool. Denote $x\stackrel{r}{\leftrightarrow}y$ as a short for $d_{\chem
}(x,y)\le r$.
Denote by $B_{\chem}(r)$ the ball in the intrinsic distance, that is, the
random set defined by $B_{\chem}(r)=\{x\dvtx0\stackrel{r}{\leftrightarrow
}x\}$.
Denote $G(r)=\mathbb{E}|B_{\chem}(r)|$. This will be our main object
of study, and we will get better and better estimates for it. As before,
assume for the sake of this sketch that $d=3$.

\subsubsection*{\texorpdfstring{Step 1.}{Step 1}}

In any transitive graph, $G(r)=e^{o(r)}$. This is due to Russo's
formula, since the number of pivotal edges for the event $0\stackrel
{r}{\leftrightarrow}x$
is always $\mbox{$\le$} r$---only edges on the path can be pivotal! See
Lemma \ref{lem:Goer}, page \pageref{lem:Goer}.

\subsubsection*{\texorpdfstring{Step 2.}{Step 2}}

We apply Schramm's lemma to $T\times T$ as explained in Section~\ref{sub:What},
and get that $\mathbb{P}(0\leftrightarrow x)\le C|x|^{2}2^{-|x|/2}$.
See Lemma \ref{lem:Oded}, page \pageref{lem:Oded}.

\subsubsection*{\texorpdfstring{Step 3.}{Step 3}}

We now repeat the argument of \cite{KN}, Theorem 1.2(i). The claim
there was that under the triangle condition $G(r)\le Cr$. The skeleton
of the argument is as follows. It is enough to prove that
%
\begin{equation}\label{eq:G2rKN}
G(2r)\ge c\frac{G(r)^{2}}{r},
\end{equation}
because if $G(r)>(2/c)r$ then it starts growing exponentially, contradicting
the information we gathered at step 1. Denote for the purpose of this
sketch by $w$ a vertex of the graph quite close to 0. $w$ is the
``opening,'' a standard step in any use of the triangle condition.
Examine the following quantity
\[
\mathbb{E}\bigl|\bigl\{(x,y)\dvtx0\stackrel{r}{\leftrightarrow}x,xw\stackrel
{r}{\leftrightarrow}y,x\nleftrightarrow xw\bigr\}\bigr|,
\]
where $xw$ stands for the product in the relevant group. You should
think about $x$ as being ``roughly pivotal'' for the event
$0\leftrightarrow y$,
with the measure of roughness related to $w$. A calculation using
the Aizenman and Newman ``off method'' which may be found in
\cite{KN}, Lemma 3.2,
or in Lemma \ref{lem:AO1_1}, page~\pageref{lem:AO1_1} below, gives
%
\begin{eqnarray}\label{eq:offintro}
&&\mathbb{E}\bigl|\bigl\{(x,y)\dvtx0\stackrel{r}{\leftrightarrow}x,xw\stackrel
{r}{\leftrightarrow}y,x\nleftrightarrow
xw\bigr\}\bigr|\nonumber\\[-8pt]\\[-8pt]
&&\qquad\ge
G(r)^{2}\biggl(1-\sum_{u,v}\mathbb{P}\bigl(0\stackrel{r}{\leftrightarrow
}u\bigr)\mathbb{P}(u\leftrightarrow
v)\mathbb{P}\bigl(v\stackrel{r}{\leftrightarrow}w\bigr)\biggr).
\nonumber
\end{eqnarray}
The sum inside the parenthesis is not \textit{quite} the triangle sum
because $0\ne w$. This expression (without the $r$'s) is known as
the \textit{open triangle sum} and it is known that if the triangle
condition holds (recall that we are still in the setting of \cite{KN09}
where the triangle condition is assumed) then the open triangle sum
tends to $0$ as the point $w$ is sent\vspace*{1pt} to infinity (\cite{BA91},
Lemma 2.1 for~$\mathbb{Z}^{d}$ and \cite{K} for a general transitive
graphs).
Taking it to be $\mbox{$\le$}\frac{1}{2}$ gives
\[
\mathbb{E}\bigl|\bigl\{(x,y)\dvtx0\stackrel{r}{\leftrightarrow}x,xw\stackrel
{r}{\leftrightarrow}y,x\nleftrightarrow xw\bigr\}\bigr|\ge{ \tfrac
{1}{2}}G(r)^{2}.
\]
At this point, a simple modification argument that allows to connect
$x$ and $xw$ while paying only a constant. The modification puts
$y$ into $\mathcal{C}(0)$ and also turns $x$ from ``roughly
pivotal'' to being properly pivotal. This proves~(\ref{eq:G2rKN})---the
modification costs only a constant, but the counting over
$x$ costs another~$r$ which explains the $\nicefrac{1}{r}$ factor
in (\ref{eq:G2rKN}). This finishes the proof
in~\cite{KN09}.\looseness=-1

How does all this apply in our case? Schramm's lemma gives only that
the triangle condition grows moderately, not that it is finite. This
means simply that it is necessary to open the triangle wider. A calculation
(done in Lemma \ref{lem:openoded}, page~\pageref{lem:openoded})
shows that it is enough to take $w$ in distance $\approx\log r$.
This is good, but not as good as it sounds because once one tries
to apply\vspace*{1pt} the modification argument one loses the exponent of the distance
namely an~$r^{C}$ factor and gets instead of (\ref{eq:G2rKN})
\[
G(2r)\ge\frac{G(r)^{2}}{r^{C}},
\]
which only shows that $G(r)\le r^{C}$---a vast improvement over
$e^{o(r)}$ but still not what we want. See Lemmas \ref{lem:AO2}
and \ref{lem:AO3}, starting from page \pageref{lem:AO2}, for the
modification argument.

\subsubsection*{\texorpdfstring{Step 4.}{Step 4}}

Because of the symmetries of the graph, our newly acquired knowledge
$G(r)\le r^{C}$ allows to get much better estimates for $0\stackrel
{r}{\leftrightarrow}x$
namely because $x$ has $2^{|x|}$ clones we get
\[
\mathbb{P}\bigl(0\stackrel{r}{\leftrightarrow}x\bigr)\le r^{C}2^{-|x|},
\]
which is better than the $2^{-|x|/2}$ given by Schramm's lemma whenever
$|x|\ge C\log r$. This allows to separate 0 and $w$ by only $\approx\log
\log r$
yielding a final $G(r)\le r(\log r)^{C}$. See Lemmas \ref{lem:AO1loglog}
and \ref{lem:AO2loglog}, page \pageref{lem:AO1loglog}.

\subsubsection*{\texorpdfstring{Step 5.}{Step 5}}\label{sub:Step-5}

We now find ourselves in a rather ridiculous situation. We have a
very good estimate for the ball in the intrinsic distance---$\mathbb
{E}|B_{\chem}(r)|\le r(\log r)^{C}$---but still absolutely no estimate for the extrinsic ball, for all
we know we might have $\mathbb{E}|\{x\dvtx0\leftrightarrow x\mbox{ and
}|x|\le r\}|\simeq2^{r/2}$,
which would of course imply that it intersects an enormous intrinsic
ball, all tightly curled up. We need to contradict this possibility
and we use the property that in a~nonamenable graph, for any $x$,
if one examines the ball of radius $r$ around~$x$ then the vast
majority of it is further from $0$ than $x$. Using the fact
that~$B_{\chem}(x,r)$ is quite small and the symmetries of the graph, we
get that we can show that the parts of $B_{\chem}(x,r)$ that go ``back,''
that is, are closer to~$0$ than $x$ are dominated by a subcritical
process. We get that the process is ``ballistic'' in the sense
that points $x$ with $d_{\chem}(0,x)=r$ have also $|x|\approx r$.
This shows that
\[
\mathbb{E}\bigl|\bigl\{x\dvtx0\leftrightarrow x\mbox{ and }|x|\le r\bigr\}\bigr|\le Cr^{3}
\]
(the 3 is just an artifact of sloppiness) and by the symmetries of
the graph one last time
\[
\mathbb{P}(0\leftrightarrow x)\le C|x|^{3}2^{-|x|},
\]
which shows the triangle condition by a direct calculation. See Lemmas
\ref{lem:ext} and \ref{lem:exiclm}, starting page \pageref{lem:ext}.

\subsection{\texorpdfstring{The case of $T\times\mathbb{Z}$.}{The case of $T\times\mathbb{Z}$}}\label{sub:TtimesZ}

Where does all this break for $T\times\mathbb{Z}$? We used the existence
of clones in every step. Hence, this argument cannot give any estimate
for $\mathbb{P}(0\leftrightarrow(0,n))$ where $(0,n)$ stands for
a vertex in the same copy of $\mathbb{Z}$ as $0$. This does not
seem like a big deal because there are not many of those. But in fact,
it breaks the argument at step 4. In other words, you can get that
$G(r)\le r^{C}$ but cannot progress beyond that, which breaks the
final step, the subcriticality of the backward process.

\begin{figure}

\includegraphics{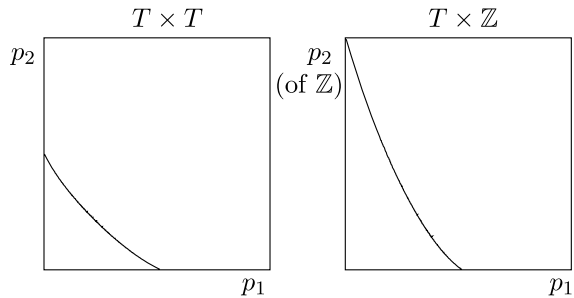}

\caption{The critical curves for $T\times T$ and $T\times\mathbb{Z}$.
Values were calculated by an invasion percolation algorithm run until
cluster size reached $10^{6}$.}\label{fig:TtimesZ}
\end{figure}

There is a different way to view this. Let our percolation have different~$p$
in the different coordinates. We get a two-parameter family of
processes with a critical curve separating the regime of only finite
clusters and the regime of infinite clusters. See Figure \ref
{fig:TtimesZ}.%
%
We immediately note the following difference between the $T\times T$
and $T\times\mathbb{Z}$ case. In the $T\times T$ case, $\nabla$
is bounded uniformly on the critical curve. On $T\times\mathbb{Z}$
it diverges as you approach the point $(0,1)$. Similarly, $\mathbb
{P}(0\leftrightarrow(0,n))$
is not bounded away from 1 uniformly on the critical curve, it converges
to $1$ as you approach $(0,1)$. I do not claim that this is a~significant
hurdle, just note that all ideas in this paper (including Schramm's
lemma) work uniformly on the entire critical line so perhaps a new
idea is needed.

\subsection{\texorpdfstring{Notation and conventions.}{Notation and conventions}}\label{sub:Notations}

All graphs in this paper will have one vertex denoted by 0. For Cayley
graphs, 0 will be the identity element. We denote by $d(x,y)$ the
graph distance of $x$ from $y$, that is, the the number of edges in
the shortest path between $x$ and $y$. Denote also $|x|=d(x,0)$
and balls by $B(x,r)=\{y\dvtx d(x,y)\le r\}$ and $B(r)=B(0,r)$. For a
subset of vertices $A$, we denote by $\partial A$ the set of all
edges with one vertex in $A$ and one outside $A$.

For percolation, we denote by $d_{\chem}(x,y)$ the length of the shortest
open path between $x$ and $y$, or $\infty$ if $x\nleftrightarrow y$.
We denote by $x\stackrel{r}{\leftrightarrow}y$ the event $d_{\chem
}(x,y)\le r$.
We denote $B_{\chem}(x,r)=\{y\dvtx x\stackrel{r}{\leftrightarrow}y\}$
and $B_{\chem}(r)=B_{\chem}(0,r)$. Be careful not to confuse $B(r)$
(which is a deterministic quantity) and $B_{\chem}(r)$ which is a
random variable. We denote $G(r)=\mathbb{E}|B_{\chem}(r)|$. Denote
also the triangle sum
\[
\nabla=\sum_{u,v}\mathbb{P}(0\leftrightarrow u)\mathbb
{P}(u\leftrightarrow v)\mathbb{P}(0\leftrightarrow v)
\]
and the restricted open triangle sum with opening $w$ and distance
$r$
\[
\nabla(w;r)=\sum_{u,v}\mathbb{P}\bigl(0\stackrel{r}{\leftrightarrow}u\bigr)\mathbb
{P}(u\leftrightarrow v)\mathbb{P}\bigl(v\stackrel{r}{\leftrightarrow}w\bigr).
\]
More standard percolation notations used are $\mathcal{C}(x)=\{
y\dvtx y\leftrightarrow x\}$,
and $A\circ B$ for the event that $A$ and $B$ ``occur disjointly,''
see \cite{G99}, Section 2.3, for the notation and for the van den Berg--Kesten
inequality $\mathbb{P}(A\circ B)\le\mathbb{P}(A)\mathbb{P}(B)$. We
shall denote the van den Berg--Kesten inequality by BK for short. Denote
Harris' inequality (\cite{G99}, Section 2.2) by FKG. For all these notations,
unless~$p$ is specified explicitly it is taken to be the $p_{c}$
of the relevant graph. When we want to examine a different $p$ we
will use the notations $\mathbb{P}_{p}$ and $\mathbb{E}_{p}$ for
the probability and the expectation with respect to $p$.

We denote by $T$ a regular tree of degree $d\ge3$. For a vertex
$x\in T\times T$, denote by $x_{1}$ and $x_{2}$ its two coordinates.
It is easy to verify that $|x|=|x_{1}|+|x_{2}|$. We denote by $\Gamma$
the group whose Cayley graph is $T\times T$ (see Lemma \ref{lem:BLPS}).

Bold letters will be used for the high-degree tree $\mathbf{T}$ that
appears in the proof of Schra\-mm's lemma (both the sketch in Section
\ref{sub:Oded}
above and the proof of Lemmas~\ref{lem:invariance} and \ref{lem:Oded}
below) and for vertices and subgraphs of it. Be careful not to
confuse~$\mathbf{T}$, which is a tree of degree $m$ and is just an auxiliary
object, with $T$ which is a tree of degree $d$ and the principle
object of investigation.

By $c$ and $C$, we will denote constants which depend only on our
graph~$G$ (usually this is $T\times T$ so they only depend on the
degree of $T$, but in Section~\ref{sec:gene} the results are general).
$c$ will denote constants which are ``small enough'' and $C$
constants which are ``big enough.'' $C$ and $c$ may refer to
different constants in different formulas and even within the same
formula. We will sometimes number them for clarity. A notation like
$C_{187}$ is specific to the lemma in which it appears. When a probability
decays exponentially in some parameter $n$, will usually denote it
by \mbox{$<$}$2e^{-cn}$.

The notation $A\approx B$ means that some constants $c$ and $C$
exist such that $cB\le A\le CB$. The notations $\simeq$ and $\lesssim$
mean nothing in particular. We only use them when we want to indicate
that two quantities are heuristically similar, but do not want to
indicate in which sense exactly. The notation $X\sim Y$ for two random
variables means that they have the same distribution. For a real number
$x$, $\lceil x\rceil$ will denote the smallest integer
$\ge x$.

\section{\texorpdfstring{Preliminaries for transitive graphs.}{Preliminaries for transitive graphs}}\label{sec:gene}
\begin{theorem}[(Aizenman and Barsky)]\label{thm:AB}For any vertex-transitive graph
$G$ and any $p<p_{c}(G)$, we have $\mathbb{E}_{p}|\mathcal{C}(0)|<\infty$.
\end{theorem}

Aizenman and Barsky \cite{AB87} formulated their result only for $\mathbb
{Z}^{d}$,
but it is well known that it holds for any transitive graph. For example,
it is mentioned in passing in \cite{S01}. A proof may be found in
\cite{AV08} or in \cite{Kb}, Appendix A.
\begin{lem}
\label{lem:Goer}For any vertex transitive graph $G$, $G(r)=e^{o(r)}$.
\end{lem}

[Recall that $G(r)=\mathbb{E}_{p_{c}}|B_{\chem}(r)|$.]
\begin{pf*}{Proof of Lemma \ref{lem:Goer}}
Fix some $x\in B(r)$ and examine the event $0\stackrel
{r}{\leftrightarrow}x$.
By Russo's formula (\cite{G99}, Section 2.4), for any $0<p<1$,
\[
\frac{d}{dp}\mathbb{P}_{p}\bigl(0\stackrel{r}{\leftrightarrow}x\bigr)=\frac
{1}{p}\mathbb{E}_{p}(|\{\mbox{open pivotal edges}\}|).
\]
We are allowed to use Russo's formula, since this event is determined
by a finite number edges, namely those of $B(r)$. For any configuration
where $0\stackrel{r}{\leftrightarrow}x$, the number of pivotal edges
is $\mbox{$\le$}r$ since clearly any edge off the path between $0$ and $x$
is not pivotal. Hence, we get
\[
\frac{d}{dp}\mathbb{P}_{p}\bigl(0\stackrel{r}{\leftrightarrow}x\bigr)\le\frac
{r}{p}\mathbb{P}_{p}\bigl(0\stackrel{r}{\leftrightarrow}x\bigr).
\]
Summing over $x\in B(r)$, we get
\[
\frac{d}{dp}\mathbb{E}_{p}|B_{\chem}(r)|\le\frac{r}{p}\mathbb
{E}_{p}|B_{\chem}(r)|
\]
or
\[
\frac{d}{dp}\log\mathbb{E}_{p}|B_{\chem}(r)|\le\frac{r}{p}.
\]
Assume by contradiction that $G(r)\ge e^{cr}$ for some $c>0$ and
infinitely ma\-ny~$r$'s. We get for $p<p_{c}$
\[
\log\mathbb{E}_{p}|B_{\chem}(r)|\ge cr-\frac{r}{p}(p_{c}-p)
\]
so for $p\in(p_{c}/(1+c),p_{c})$ we have $\mathbb{E}_{p}|B_{\chem
}(r)|\to\infty$
as $r\to\infty$. This contradicts the theorem of Aizenman and Barsky.
\end{pf*}
\begin{lem}
\label{lem:diagram}For any transitive graph, any $r>0$ and any $\lambda>0$,
\[
\mathbb{P}\bigl(|B_{\chem}(r)|>\lambda G(r)^{2}\bigr)\le2e^{-c\lambda},
\]
where $c$ is an absolute constant.
\end{lem}
\begin{pf}
This is a standard corollary of Aizenman and Newman's diagrammatic bounds.
Refer to \cite{G99}, Section 6.3, for a complete treatment. The classic
picture is as follows. We wish to calculate the $n$th
moment of $|B_{\chem}(r)|$. For this, we note that if $0\leftrightarrow x_{i}$
for $i=1,\ldots,n$ then there exist $y_{1},\ldots,y_{n-1}$ and a
tree describing the connection scheme $U$ with $n+1$ leaves (corresponding
to the vertices $0,x_{1},\ldots,x_{n}$) and $n-1$ inner points (corresponding
to the vertices $y_{1},\ldots,y_{n-1}$) all of which are of degree
3 such that for every edge of $U$, the corresponding vertices are
connected by an open path, and all such paths are edge-disjoint. For
convenience, denote $(z_{1},\ldots,z_{2n})=(0,x_{1},\ldots
,x_{n},y_{1},\ldots,y_{n-1})$
and let the vertices of the tree $U$ be the numbers $\{1,\ldots,2n\}$
correspondingly.

Now, in our case the paths from $0$ to $x_{i}$ are constrained to
be of length $\le r$ and this constraint is carried over to all internal
paths of our tree. Hence,
\begin{eqnarray*}
&&\phantom{\mbox{By BK}}\mathbb{E}\bigl|\bigl\{(x_{1},\ldots,x_{n})\dvtx0\stackrel{r}{\leftrightarrow}x_{i}
\ \forall
i\bigr\}\bigr|  \\
&&\phantom{\mbox{By BK}}\qquad\le\sum_{U}\sum_{z_{2},\ldots,z_{2n}}\mathbb{P}\bigl(z_{i}\stackrel
{r}{\leftrightarrow}z_{j}\mbox{ for all edges }(i,j)\mbox{ of }U\mbox{,
disjointly}\bigr)\\
&&\mbox{By BK}\qquad \le\sum_{U}\sum_{z_{2},\ldots,z_{2n}}\prod
_{(i,j)\in U}\mathbb{P}\bigl(z_{i}\stackrel{r}{\leftrightarrow}z_{j}\bigr)\\
&&\phantom{\mbox{By BK}}\qquad=\sum_{U}G(r)^{2n-1}=1\cdot3\cdot5\dotsb(2n-3)\cdot G(r)^{2n-1}\\
&&\phantom{\mbox{By BK}}\qquad\le(2nG(r)^{2})^{n}.
\end{eqnarray*}
The proof is now complete as
\begin{eqnarray*}
\mathbb{P}\bigl(|B_{\chem}(r)|>\lambda
G(r)^{2}\bigr)&=&\mathbb{P}\bigl(|B_{\chem}(r)|^{n}>\lambda^{n}G(r)^{2n}\bigr)\\
&\le&\frac{\mathbb{E}(|B_{\chem}(r)|^{n})}{\lambda^{n}G(r)^{2n}}\le
\biggl(\frac{2n}{\lambda}\biggr)^{n}
\end{eqnarray*}
and setting $n=\lfloor\lambda/4\rfloor$ we are done.
\end{pf}

The last general claim we wish to demonstrate before moving on to
the product of two trees is the invariance step in Schramm's lemma.
Recall the sketch of the lemma and also the discussion in Section \ref
{sub:What}.
The claim there is that for any branching process with a time-reversible
step, the resulting~$\mathbf{W}$ is invariant to the automorphisms
of the tree $\mathbf{T}$. We will now prove this fact, and we formulate
it in sufficient generality so that it can be used both for the original
version of Schramm's lemma and for our purposes.

This is the only place in the paper where it matters which Cayley
graph we are talking about, so let us fix that all Cayley graphs are
\textit{right} Cayley graphs, that is, if $\Gamma$ is a finitely generated
group and $S$ a set of generators then the edges of the Cayley graphs
are $\{(g,gs)\dvtx g\in\Gamma,s\in S\}$.
\begin{defn}
\label{def:oded}Let $\Gamma$ be a finitely generated group and let
$S$ be a set of generators. Let $\mu$ be some discrete measure on
$\Gamma$ with $\mu(x)=\mu(x^{-1})$. Let $m\ge1$ be some integer.
Let $0<p<1$.
\begin{itemize}
\item Define $\mathbf{T}$ to be a regular tree of degree $m+1$. We will
use $\mathbf{T}$ also to denote the set of vertices of $\mathbf{T}$,
and the edges will be denoted by $E(\mathbf{T})$. Fix one element
of~$\mathbf{T}$, call it the root and denote it by $\mathbf{r}$.
\item Define $\pi\dvtx\mathbf{T}\to\Gamma$ which is a random map (``the
locations of the particles''). For $\mathbf{r}$ the root of the
tree $\mathbf{T}$, we define $\pi(\mathbf{r})=0$ where $0$ is the
identity element of~$\Gamma$. We continue inductively. Assume $\pi
(\mathbf{v})$
is already defined. For every child $\mathbf{w}$ of $\mathbf{v}$,
we define $\pi(\mathbf{w})=\pi(\mathbf{v})X_{\mathbf{v},\mathbf{w}}$
where $X_{\mathbf{v},\mathbf{w}}$ are i.i.d. random variables distributed
like $\mu$.
\item Finally, define Schramm's process $\mathbf{W}=\mathbf{W}(\Gamma
,S,\mu,m,p)$
to be a random subset of the edges of $\mathbf{T}$ defined by
\[
(\mathbf{u},\mathbf{v})\in\mathbf{W}\quad\iff\quad\pi(\mathbf{u})\leftrightarrow
\pi(\mathbf{v}) \qquad\forall(\mathbf{u},\mathbf{v})\in
E(\mathbf{T}),
\]
where $g\leftrightarrow h$ denotes that $g$ and $h$ are connected
in a $p$-percolation process independent of the $X_{\mathbf{u},\mathbf{v}}$
on the (right) Cayley graph of $\Gamma$ with respect to the set of
generators $S$.
\end{itemize}
\end{defn}
\begin{lem}
\label{lem:invariance}For any $\Gamma$, $S$, $\mu$, $m$ and $p$
as above, Schramm's process is invariant to the automorphisms of
$\mathbf{T}$.
\end{lem}
\begin{pf}
Denote the root of $\mathbf{T}$ by $\mathbf{r}$ and its children
by $\mathbf{v}_{1},\ldots,\mathbf{v}_{m+1}$. Let $\varphi\dvtx\mathbf{T}\to
\mathbf{T}$
be the automorphism that takes $\mathbf{v}_{1}\mapsto\mathbf{r}$
and $\mathbf{r}\mapsto\mathbf{v}_{m+1}$ but otherwise preserves the
order among children: the $m$ children of $\mathbf{v}_{1}$ are mapped
to $\mathbf{v}_{1},\ldots,\mathbf{v}_{m}$ in order, $\mathbf
{v}_{2},\ldots,\mathbf{v}_{m+1}$
are mapped to the $m$ children of $\mathbf{v}_{m+1}$ in order etc.
For every permutation of $m+1$ elements $\sigma$, let $\psi_{\sigma
}\dvtx\mathbf{T}\to\mathbf{T}$
be the automorphism permuting the children of $\mathbf{r}$ according
to $\sigma$ but otherwise preserving the order. Let $H$ be the group
of automorphisms of $\mathbf{T}$ generated by $\varphi$ and all
$\psi_{\sigma}$.

It is straightforward to verify that for any automorphism $\theta$
of $\mathbf{T}$ and every~$r$ there exists an $\eta\in H$ which
is identical to $\theta$ on the entire ball of radius $r$ in $\mathbf{T}$.
In other words, $H$ is dense in the compact-open topology. This implies
that it is enough to show that the distribution of $\mathbf{W}$ is
invariant to the action of $H$ and hence it is enough to show that
it is invariant to the action of~$\varphi$ and~$\psi_{\sigma}$.
Verifying $\psi_{\sigma}$ is immediate, so we are left with showing
that $\varphi\mathbf{W}\sim\mathbf{W}$ where we define $(\varphi\mathbf
{W})(e)=\mathbf{W}(\varphi^{-1}(e))$.
It will be more convenient to verify that $\varphi^{-1}\mathbf{W}\sim
\mathbf{W}$
and we will do so.

Write our probability space as $\Omega_{1}\times\Omega_{2}$ where
$\Omega_{1}$ is the probability space of the bran\-ching random
walk and $\Omega_{2}$ is the probability space of the percolation.
Further, write $\Omega_{1}$ as $\Gamma^{E(\mathbf{T})}$ where $E(\mathbf{T})$
is the set of edges of $\mathbf{T}$, with the measure being the product
measure $\mu^{E(\mathbf{T})}$. Any automorphism $\varphi$ of $\mathbf{T}$
induces a measure preserving map $\alpha\dvtx\Omega_{1}\to\Omega_{1}$
by
\[
\alpha(\omega)(\mathbf{x},\mathbf{y})=\cases{
\omega(\varphi(\mathbf{x}),\varphi(\mathbf{y})), &\quad $\varphi$
preserves the orientation of
$(\mathbf{x},\mathbf{y})$,\cr
\omega(\varphi(\mathbf{y}),\varphi(\mathbf{x}))^{-1}, &\quad otherwise.}
\]
The $^{-1}$ on the bottom clause stands for inversion in the group
$\Gamma$. Also we need to explain what does it mean that ``$\varphi$
preserves the orientation of $(\mathbf{x},\mathbf{y})$''---this
means, when $\mathbf{x}$ is the father of $\mathbf{y}$,
that $\varphi(\mathbf{x})$ is the father of~$\varphi(\mathbf{y})$.
Of course, our $\varphi$ only reverses the orientation of one edge,
$(\mathbf{r},\mathbf{v}_{1})$. $\alpha$~is measure preserving
because the measure $\mu$ is invariant to the operation~$^{-1}$.
Slightly abusing notations we consider $\alpha$ also as a map $\Omega
_{1}\times\Omega_{2}\to\Omega_{1}\times\Omega_{2}$
acting only on the first coordinate.

We now define a second measure preserving map $\beta\dvtx\Omega_{1}\times
\Omega_{2}\to\Omega_{1}\times\Omega_{2}$
as follows: for any $\omega\in\Omega_{1}$ we let $f(\omega)$ be
an automorphism of the Cayley graph of $\Gamma$ given by
\[
f(\omega)v=\omega(\mathbf{r},\mathbf{v}_{m+1})^{-1}v,
\]
where the product is in the group $\Gamma$---again we consider
$\omega$ as an element of $\Gamma^{E(\mathbf{T})}$ so $\omega(\mathbf
{r},\mathbf{v}_{m+1})$
is simply the position of the $(m+1)$st child of the
original particle. We also consider $f$ as acting on $\Omega_{2}$
(which is just the product space $\{0,1\}^{E(\Gamma)}$) by $f\omega
_{2}(v,w)=\omega_{2}(f^{-1}(v),f^{-1}(w))$.
We now define $\beta(\omega_{1},\omega_{2})=(\omega_{1},f(\omega
_{1})\omega_{2})$.
Since $f(\omega_{1})$ is measure preserving for any $\omega_{1}$,
we get that $\beta$ is measure preserving by Fubini's theorem.

The lemma is now finished because applying the measure preserving
transformation $\alpha\circ\beta$ to the probability space is the
same as applying $\varphi^{-1}$ to~$\mathbf{W}$. Let us verify this
formally. We consider $\mathbf{W}$ as a function $\Omega_{1}\times\Omega
_{2}\to\{0,1\}^{E(\mathbf{T})}$
defined by ``$\mathbf{W}(\omega_{1},\omega_{2})(\mathbf{x},\mathbf{y})=1$
if $\prod\omega_{1}(\mathbf{x}_{i},\mathbf{x}_{i+1})$ is connected
to $\prod\omega_{1}(\mathbf{y}_{i},\mathbf{y}_{i+1})$ in the configuration
$\omega_{2}$,'' where for an element $\mathbf{x}\in\mathbf{T}$
we define $\mathbf{x}_{0},\mathbf{x}_{1},\ldots$ to be the elements
of the tree on the branch from $\mathbf{r}=\mathbf{x}_{0}$ to $\mathbf{x}$,
and where the $\prod$ is in the group $\Gamma$ and is taken
left-to-right,
that is, $\omega(\mathbf{x}_{0},\mathbf{x}_{1})\omega(\mathbf
{x}_{1},\mathbf{x}_{2})\dotsb.$
We wish to show that $\mathbf{W}(\alpha(\beta(\omega)))=\varphi
^{-1}\mathbf{W}(\omega)$.
But $\mathbf{W}(\alpha(\beta(\omega))(\mathbf{x},\mathbf{y})=1$ if
$\prod\alpha(\beta(\omega))_{1}(\mathbf{x}_{i},\mathbf{x}_{i+1})$
is connected to $\prod\alpha(\beta(\omega))_{1}(\mathbf{y}_{i},\mathbf
{y}_{i+1})$
in $\alpha(\beta(\omega))_{2}$. Now,
\begin{eqnarray*}
&&\alpha(\beta(\omega))_{1}(\mathbf{x},\mathbf{y})\\
&&\qquad =\cases{
\beta(\omega)_{1}(\varphi(\mathbf{x}),\varphi(\mathbf{y})),
&\quad $\varphi$ preserves the orientation of
$(\mathbf{x},\mathbf{y})$,\cr
\beta(\omega)_{1}(\varphi(\mathbf{y}),\varphi(\mathbf{x})
)^{-1}, &\quad otherwise,}
\\
&&\qquad =\cases{
\omega_{1}(\varphi(\mathbf{x}),\varphi(\mathbf{y})), &\quad
$\varphi$ preserves the orientation of
$(\mathbf{x},\mathbf{y})$\cr
\omega_{1}(\varphi(\mathbf{y}),\varphi(\mathbf{x}))^{-1}, &\quad
otherwise,}
\end{eqnarray*}
so
\begin{eqnarray*}
&&\prod\alpha(\beta(\omega))_{1}(\mathbf{x}_{i},\mathbf{x}_{i+1})\\
&&\qquad=\cases{
\displaystyle \omega_{1}(\mathbf{r},\mathbf{v}_{m+1})^{-1}\prod_{i=1}\omega_{1}
(\varphi(\mathbf{x}_{i}),\varphi(\mathbf{x}_{i+1})),
&\quad $\mathbf{x}_{1}=\mathbf{v}_{1}$,\cr
\displaystyle \prod_{i=0}\omega_{1}(\varphi(\mathbf{x}_{i}),\varphi(\mathbf
{x}_{i+1})), &\quad otherwise.}
\end{eqnarray*}
Comparing to the branch from $\mathbf{r}$ to $\varphi(\mathbf{x})$,
we get
\[
\prod\alpha(\beta(\omega))_{1}(\mathbf{x}_{i},\mathbf{x}_{i+1})=\omega
_{1}(\mathbf{r},\mathbf{v}_{m+1})^{-1}\prod\omega_{1}(\varphi
(\mathbf{x})_{i},\varphi(\mathbf{x})_{i+1})
\]
in both cases.\vadjust{\eject}

For the percolation configuration we have a similar calculation,
\begin{eqnarray*}
\alpha(\beta(\omega))_{2}(v,w)
&=&\beta(\omega)_{2}(v,w)=(f(\omega_{1})\omega_{2})(v,w)\\
&=& \omega_{2}(f(\omega_{1})^{-1}v,f(\omega_{1})^{-1}w)\\
&=& \omega_{2}(\omega_{1}(\mathbf{r},\mathbf{v}_{m+1})v,\omega
_{1}(\mathbf{r},\mathbf{v}_{m+1})w).
\end{eqnarray*}
We put the formulas for $\Omega_{1}$ and $\Omega_{2}$ together and get that
$\mathbf{W}(\alpha(\beta(\omega)))(\mathbf{x},\allowbreak\mathbf{y})=1$ if and only if
$\omega_{1}(\mathbf{r},\mathbf{v}_{m+1})^{-1}\prod\omega_{1}(\varphi(\mathbf{x})_{i},\varphi(\mathbf{x})_{i+1})$
and $\omega_{1}(\mathbf{r},\mathbf{v}_{m+1})^{-1}\allowbreak\prod\omega_{1}(\varphi(\mathbf{y})_{i}$,
$\varphi(\mathbf
{y})_{i+1})$
are connected in the configuration $\omega_{2}(\omega_{1}(\mathbf
{r},\mathbf{v}_{m+1})v,\allowbreak\omega_{1}(\mathbf{r},\mathbf{v}_{m+1})w)$.
The terms $\omega_{1}(\mathbf{r},\mathbf{v}_{m+1})$ now cancel (recall
the a right Cayley graph is invariant to left translations) and we
get that this happens if and only if $\prod\omega_{1}(\varphi
(\mathbf{x})_{i},\varphi(\mathbf{x})_{i+1})$
is connected to $\prod\omega_{1}(\varphi(\mathbf{y})_{i},\varphi
(\mathbf{y})_{i+1})$
in $\omega_{2}$ which is exactly $\mathbf{W}(\omega)(\varphi
(\mathbf{x}),\varphi(\mathbf{y}))$.
In other words
\[
\mathbf{W}(\alpha(\beta(\omega))(\mathbf{x},\mathbf{y})=\mathbf
{W}(\omega)(\varphi(\mathbf{x}),\varphi(\mathbf{y}))
\qquad\mbox{for all }\mathbf{x}\mbox{ and }\mathbf{y},
\]
which is exactly $\mathbf{W}(\alpha(\beta(\omega)))=\varphi^{-1}\mathbf
{W}(\omega)$
which shows that the distribution of $\mathbf{W}$ is invariant to
$\varphi^{-1}$ and hence to $\varphi$. As explained, this shows that
$\mathbf{W}$ is invariant to a group of automorphisms dense in the
compact-open topology, hence to all automorphisms, proving the lemma.
\end{pf}

\section{\texorpdfstring{The product of two trees.}{The product of two trees}}
\begin{lem}
\label{lem:BLPS}At $p_{c}$ there is no infinite cluster for $T\times T$.
\end{lem}
\begin{pf}
By \cite{BLPS99}, every nonamenable Cayley graph satisfies this property.
Hence, we need only show that $T\times T$ is a nonamenable Cayley
graph. This however is easy. Any tree of degree $d$ is the Cayley
graph of the free product
\[
G_{d}:=\underbrace{\nicefrac{\mathbb{Z}}{2\mathbb{Z}}*\dotsb*\nicefrac
{\mathbb{Z}}{2\mathbb{Z}}}_{d\ \mathrm{times}}
\]
with the natural generators (namely one from every copy of $\nicefrac
{\mathbb{Z}}{2\mathbb{Z}}$).
Of course if $d$ is even then one can simply take a free group with
$\frac{1}{2}d$ generators. The product is thus the Cayley graph of
the group $G_{d}\times G_{d}$. The claim of nonamenability is just
as easy. We follow \cite{BLPS99} and say that a graph is nonamenable
if there exists some $c$ such that for all finite $A$, $|\partial A|\ge c|A|$.
Let therefore $A\subset T\times T$ be any finite set. For any $x\in T$,
denote by $A_{x}$ the slice $\{y\dvtx(y,x)\in A\}$. Then, when $A_{x}\ne
\varnothing$,
%
\begin{equation}\label{eq:treeisop}
|\partial A_{x}|=(d-2)|A_{x}|+2,
\end{equation}
where $\partial$ is the edge boundary in the tree $T$. Equation (\ref{eq:treeisop})
is a well-known property of regular trees and may be readily proved
by induction on $|A_{x}|$. Summing over $x$, we get
\[
|\partial A|\ge\sum_{x}|\partial A_{x}|>\sum_{x}(d-2)|A_{x}|=(d-2)|A|
\]
as needed.\vadjust{\eject}
\end{pf}
\begin{defn}
We denote by $\Gamma$ the group whose Cayley graph is $T\times T$,
namely $G_{d}\times G_{d}$ from the previous lemma.
\end{defn}

The next lemma is the adaptation of Schramm's argument to our
setting.\vspace*{6pt}
\begin{lem}
\label{lem:Oded}For any $x\in T\times T$,
\[
\mathbb{P}(0\leftrightarrow x)\le C|x|^{2}(d-1)^{-|x|/2}.
\]
\end{lem}
\begin{pf}
Let $x=(x_{1},x_{2})$ and let $k_{i}=|x_{i}|$ be the distances of
the coordinates of $x$ from the root in the two trees. For any
$y=(y_{1},y_{2})\in T\times T$
denote
\[
L(y)=\{(z_{1},z_{2})\dvtx d(z_{1},y_{1})=k_{1},d(z_{2},y_{2})=k_{2}\}.
\]
Clearly $|L(y)|=(d-1)^{|x|}$. Recall the definition of Schramm's
process (Definition~\ref{def:oded} on page \pageref{def:oded}). First,
we need a group and we take the group to be our~$\Gamma$, the group
whose Cayley graph is $T\times T$. The next element in Schramm's
process is a branching process with a time-reversible step $\mu$.
We take
\[
\mu(y)=\cases{
\dfrac{1}{|L(0)|}, &\quad $y\in L(0)$,\vspace*{2pt}\cr
0, &\quad otherwise.}
\]
Clearly, $\mu(y)=\mu(y^{-1})$, that is, is time-reversible. Finally,
we need two parameters, $m$ the branching number (which we leave
unspecified for a while) and $p$, which we take to be $p_{c}(T\times T)$.
Schramm's process is now a random subset $\mathbf{W}$ of the edges
of a regular tree of degree $m+1$ which we denote by~$\mathbf{T}$.
The statement of Lemma \ref{lem:invariance} now says
\[
\mathbf{W}\mbox{ is invariant to the automorphisms of }\mathbf{T}.
\]
Examine now the parameter $m$---recall that Schramm's process
involves a~branching random walk on $T\times T$ where each particle
splits into $m$ children and then each child makes one step of $\mu$.
For every $m$ for which the branching process is transient (i.e., with
probability $1$ only a finite number of particles return to any given
point), then the configuration $\mathbf{W}$ contains no infinite
component. This is due to Lemma \ref{lem:BLPS}, since the cluster
of $0$ is finite (with probability $1$) and only a finite number
of particles return to each of its points (again with probability
$1$). Let $\mathbf{W}$ be a subset of edges of $\mathbf{T}$, and
denote by $\mathcal{C}(\mathbf{x})=\mathcal{C}(\mathbf{x};\mathbf{W})$
the cluster of some $\mathbf{x}\in\mathbf{T}$ in $\mathbf{W}$, that is, all
$\mathbf{y}$ connected to $\mathbf{x}$ by a path of edges in $\mathbf{W}$.
Denote also $\deg\mathbf{x}=|\{\mathbf{y}\dvtx(\mathbf{x},\mathbf{y})\in
\mathbf{W}\}|$,
the degree of $\mathbf{x}$ in $\mathbf{W}$. We now define
\[
M(\mathbf{x},\mathbf{y};\mathbf{W})=\cases{
\dfrac{\deg\mathbf{x}}{|\mathcal{C}(\mathbf{x})|}, &\quad $\mathbf{y}\in
\mathcal{C}(\mathbf{x})$ and
$|\mathcal{C}(\mathbf{x})|<\infty$,\vspace*{2pt}\cr
0, &\quad otherwise.}
\]
Clearly, $M$ is invariant to the automorphisms of $\mathbf{T}$ in the sense that
$M(\mathbf{x},\mathbf{y};\allowbreak\mathbf{W})=M(\varphi\mathbf{x},\varphi\mathbf{y};\varphi\mathbf{W})$
for any automorphism $\varphi$ of $\mathbf{T}$. A simple change
of variables known as the ``mass transport principle,'' see \cite{BLPS99},
equation (2.1), shows that
\[
\sum_{\mathbf{y}\in\mathbf{T}}\mathbb{E}M(\mathbf{x},\mathbf{y};\mathbf
{W})=\sum_{\mathbf{y}\in\mathbf{T}}\mathbb{E}M(\mathbf{y},\mathbf
{x};\mathbf{W}) \qquad\forall\mathbf{x}\in\mathbf{T}
\]
($\mathbb{E}$ here is with respect to $\mathbf{W}$ and we use the
invariance of $\mathbf{W}$ too). Now, the left-hand side is obviously
just $\mathbb{E}(\deg\mathbf{x}\cdot\mathbf{1}_{|\mathcal{C}(\mathbf
{x})|<\infty})$
(the sum and the expectation may be exchanged since $M$ is positive).
Since our clusters are finite a.s., the left-hand side is simply
$\mathbb{E}(\deg\mathbf{x})$.
The right-hand side, on the other hand, is the expected average degree
of $\mathcal{C}(\mathbf{x})$. However, $\mathbf{T}$ is a~tree, and
therefore $\mathcal{C}(\mathbf{x})$ is a finite tree and it is well
known (and easy to prove inductively) that the average degree of any
finite tree is $\mbox{$<$}2$. Hence, we get\looseness=1
\[
\mathbb{E}(\deg\mathbf{x})<2.\looseness=0
\]
What is the meaning of $\mathbb{E}\deg\mathbf{x}$? Since our branching
process consists of sending $m+1$ particles to random points in $L(0)$,
and all particles are identical, it means that for a random point
$y\in L(0)$
%
\begin{equation}\label{eq:2/m}
\mathbb{P}(0\leftrightarrow y)\le\frac{2}{m+1}.
\end{equation}
However, the symmetries of the two trees show that for the $y\in L(0)$,
the probabilities $\mathbb{P}(0\leftrightarrow y)$ are all equal.
Hence, we get (\ref{eq:2/m}) for \textit{any} $y\in L(0)$. This is
the crux of Schramm's argument, and we need only to calculate an $m$
for which the branching process is still transient (ideally one might
want to find the maximal such $m$, but the following argument is
not precise).\looseness=-1
\begin{claim*}
With the definitions above, for $m=\frac{c}{|x|^{2}}(d-1)^{|x|/2}$
the branching process is transient.\vspace*{-3pt}
\end{claim*}
\begin{pf}
Examine some $y=(y_{1},y_{2})$. Let $z=(z_{1},z_{2})$
be a random point in~$L(y)$. We wish to understand the distribution
of $|z_{i}|$. Clearly, $z_{1}$ and $z_{2}$ are independent. Let
us therefore examine $z_{1}$. Going to distance $k_{1}$ from~$y_{1}$
is equivalent to making a nonbacktracking walk of distance $k_{1}$.
With probability $\frac{1}{d}$, the first step is in the direction
of $0$. If this happens, then at each step we have probability $\frac{1}{d-1}$
to step\vspace*{1pt} in the direction of $0$. Once we did one step away from $0$,
we can never go back. Therefore, the number of steps taken in the direction
of $0$ is dominated by a geometric variable with expectation $\frac{1}{d-2}$.
We get that for a given point $\mathbf{u}\in\mathbf{T}$ with $|\mathbf{u}|=l$,
the probability that the corresponding particle is in $0$ can be
bounded by
\begin{eqnarray*}
\mathbb{P}\bigl(\varphi(\mathbf{u})_{1}=0\bigr) &\le& \mathbb{P}\Biggl(\sum
_{i=1}^{l}\operatorname{Geom}_{i}\ge\frac{lk_{1}}{2}\Biggr)=\mathbb{P}\Biggl(\sum
_{i=1}^{\lceil lk_{1}/2\rceil+l-1}\operatorname{Bin}_{i}<l\Biggr)
\\
&\le & (d-1)^{-\lceil lk_{1}/2\rceil}\sum_{i=0}^{l-1}
\pmatrix{\lceil lk_{1}/2\rceil+l-1\cr i}\le
\bigl(Ck_{1}(d-1)^{-k_{1}/2}\bigr)^{l},
\end{eqnarray*}
where $\operatorname{Geom}_{i}$ are independent geometric variables with
expectation $\frac{1}{d-2}$ and $\operatorname{Bin}_{i}$ are independent
Bernoulli trials which give 0 with probability $\frac{1}{d-1}$
and~1 with probability $1-\frac{1}{d-1}$. The same calculation for the
other tree gives
\[
\mathbb{P}\bigl(\varphi(\mathbf{u})_{2}=0\bigr)\le
\bigl(Ck_{2}(d-1)^{-k_{2}/2}\bigr)^{l}
\]
and since the two trees are independent we get
\[
\mathbb{P}\bigl(\varphi(\mathbf{u})=0\bigr)\le
\bigl(Ck_{1}k_{2}(d-1)^{-|x|/2}\bigr)^{l}\le
\bigl(C_{1}|x|^{2}(d-1)^{-|x|/2}\bigr)^{l}.
\]
Taking $m$ with $m\cdot C_{1}|x|^{2}(d-1)^{-|x|/2}<1$, we see that
\[
\mathbb{E}\bigl(|\{\mathbf{u}\in\mathbf{T}\mbox{ s.t. }\varphi(\mathbf
{u})=0\}|\bigr)\le\sum_{l=0}^{\infty}\bigl(m\cdot
C_{1}|x|^{2}(d-1)^{-|x|/2}\bigr)^{l}<\infty.
\]
Thus, the process is transient, proving the claim. With (\ref{eq:2/m}),
this also finishes the proof of Lemma \ref{lem:Oded}.
\end{pf}
\noqed\end{pf}

As explained in the \hyperref[sec1]{Introduction}, Lemma \ref{lem:Oded} shows that
the triangle sum grows only polynomially, and the next step is to
show an ``open triangle condition'' [recall~(\ref{eq:offintro}),
page \pageref{eq:offintro}] with ``logarithmic opening.'' Here
is the precise formulation.
\begin{lem}
\label{lem:openoded}There exists some $C_{1}$ such that for any
$r\ge2$ and any $w\in T\times T$ with $|w|>C_{1}\log r$, $\nabla(w;r)\le
\frac{1}{2}$.
\end{lem}

Recall that $\nabla(w;r)=\sum_{u,v}\mathbb{P}(0\stackrel
{r}{\leftrightarrow}u)\mathbb{P}(u\leftrightarrow v)\mathbb
{P}(v\stackrel{r}{\leftrightarrow}w)$.
\begin{pf*}{Proof of Lemma \ref{lem:openoded}}
It will be convenient to replace the restriction that the path between
$0$ and $u$ is of length $\le r$ with simply $|u_{i}|\le r$ (as
before, $u_{1}$ and $u_{2}$ are the two coordinates and $|u_{i}|$
is their distance from the root of $T$), and similarly $|v_{i}|\le r+|w|$.
Denote $s=r+|w|$. Applying Lemma~\ref{lem:Oded}, we get
%
\begin{eqnarray}
\label{eq:triangvr}
\nabla(w;r) &=&\sum_{u,v}\mathbb{P}\bigl(0\stackrel
{r}{\leftrightarrow}u\bigr)\mathbb{P}(u\leftrightarrow v)\mathbb
{P}\bigl(v\stackrel{r}{\leftrightarrow}w\bigr)\nonumber\\
&\le& C\sum_{|u|,|v|\le
s}|u|^{2}(d-1)^{-|u|/2}|u-v|^{2}(d-1)^{-d(u,v)/2} \nonumber\\[-8pt]\\[-8pt]
&&\hspace*{40.8pt}{}\times|v-w|^{2}(d-1)^{-d(v,w)/2}
\nonumber\\
&\le& Cs^{6}\sum_{|u_{1}|,|u_{2}|,|v_{1}|,|v_{2}|\le
s}(d-1)^{-(|u|+d(u,v)+d(v,w))/2}\nonumber
\end{eqnarray}
and at this point the sum decomposes into a product
of one term for each tree. Totally we get
\[
\nabla(w;r)= Cs^{6}\biggl(\sum_{|u_{1}|,|v_{1}|\le
s}(d-1)^{-(|u_{1}|+d(u_{1},v_{1})+d(v_{1},w_{1}))/2}
\biggr)^{2}.
\]
Estimating the terms in the last expression is straightforward. The
subtree generated by $0$, $u_{1}$, $v_{1}$ and $w_{1}$ may take
one of the $3$ shapes in Figure \ref{fig:3diag}.%
%
\begin{figure}

\includegraphics{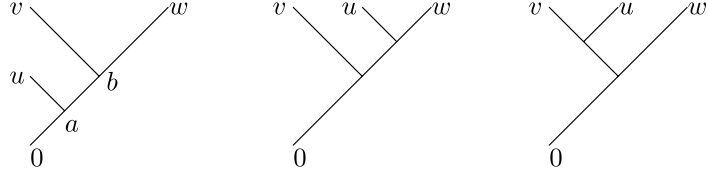}

\caption{3 ways $0$, $u$, $v$ and $w$ may be connected
in the tree.}
\label{fig:3diag}
\end{figure}
Take the first (leftmost) case as an example. Denote the two branch
points by $a$ and $b$ so that $|u_{1}|=|a|+d(u_{1},a)$,
$d(u_{1},v_{1})=d(u_{1},a)+d(a,b)+d(b,v_{1})$
and $d(v_{1},w_{1})=d(v_{1},b)+d(b,w_{1})$. We get
\[
\tfrac
{1}{2}\bigl(|u_{1}|+d(u_{1},v_{1})+d(v_{1},w_{1})\bigr)=d(a,u_{1})+d(b,v_{1})+\tfrac
{1}{2}|w_{1}|
\]
(in the second case of Figure \ref{fig:3diag} you get a $\ge$ rather
than an $=$). Fixing $a$ and summing over all $u_{1}$ gives
\[
\mathop{\sum_{u_{1}\ \mathrm{in}\ \mathrm{the}\ \mathrm{subtree}}}_
{\mathrm{of}\ a, |u_{1}-a|\le s}
(d-1)^{-d(a,u_{1})}=s.
\]
Similarly, fixing $b$ and summing over $v_{1}$ gives another factor
of $s$. Finally, $a$~and~$b$ have $\mbox{$\le$} s$ possibilities each.
Hence, we get
\[
\mathop{\sum_{u_{1},v_{1}\ \mathrm{connected}\ \mathrm{as}}}_{\mathrm
{in}\ \mathrm{the}\ \mathrm{first}\ \mathrm{diagram}}
(d-1)^{-(|u_{1}|+d(u_{1},v_{1})+d(v_{1},w_{1}))/2}\le s^{4}(d-1)^{-|w_{1}|/2}.
\]
A similar calculation works for the other 2 diagrams and we get
\[
\sum_{|u_{1}|,|v_{1}|\le
s}(d-1)^{-(|u_{1}|+d(u_{1},v_{1})+d(v_{1},w_{1}))/2}\le3s^{4}(d-1)^{-|w_{1}|/2}.
\]
The sum over $u_{2}$, $v_{2}$ and $w_{2}$ is the same, and multiplying
we get
%
\begin{equation}\label{eq:yzr}
\sum_{|u|,|v|\le s}(d-1)^{-(|u|+d(u,v)+d(v,w))/2}\le
9s^{8}(d-1)^{-|w|/2}
\end{equation}
(we remark this partial step in the computation as it will be needed
below in Lemma \ref{lem:AO1loglog}. Note that it works for an arbitrary
$s$ and not just for $s=r+|w|$). Inserting into (\ref{eq:triangvr}),
we end up with
\[
\nabla(w;r)\le C(r+|w|)^{14}(d-1)^{-|w|/2}
\]
and the lemma is finished: with a choice of $C_{1}$ sufficiently
large we get $\nabla(w;r)\le\frac{1}{2}$.
\end{pf*}
\begin{lem}
\label{lem:AO1_1}Let $C_{1}$ be as in Lemma \ref{lem:openoded}.
Then for any $r\ge2$ and any $w$ with $|w|\ge C_{1}\log r$,
\[
\mathbb{E}\bigl(\bigl|\bigl\{(x,y)\dvtx0\stackrel{r}{\leftrightarrow}x,  xw\stackrel
{r}{\leftrightarrow}y, 0\nleftrightarrow y\bigr\}\bigr|\bigr)\ge{ \tfrac
{1}{2}}G(r)^{2},
\]
where $xw$ stands for the product in the group $\Gamma$ whose Cayley
graph is $T\times T$.
\end{lem}

Note that in the restriction $0\nleftrightarrow y$ we do not require
anything from the length of the path. The restrictions that the path
is $\mbox{$\le$}r$ apply only to the paths from $0$ to $x$ and from $xw$
to $y$.

The proof is identical to that of \cite{KN09}, Lemma 3.2, and we
include it mainly for completeness.
\begin{pf*}{Proof of Lemma \ref{lem:AO1_1}}
By multiplying\vspace*{1pt} with $x^{-1}$ from the left and then doing the change
of variables $x^{-1}y\mapsto y$, $x^{-1}\mapsto x$, we see that it
is enough to show
\[
\mathbb{E}\bigl(\bigl|\bigl\{(x,y)\dvtx0\stackrel{r}{\leftrightarrow}x,  w\stackrel
{r}{\leftrightarrow}y, 0\nleftrightarrow w\bigr\}\bigr|\bigr)\ge{ \tfrac
{1}{2}}G(r)^{2}.
\]
Fix some $x$ and $y$. Now let us condition on the cluster of 0,
$\mathcal{C}(0)$. We get
\begin{eqnarray*}
&&\mathbb{P}\bigl(0\stackrel{r}{\leftrightarrow}x,
w\stackrel{r}{\leftrightarrow}y, 0\nleftrightarrow w\bigr)\\
&&\qquad=\sum_{\mathrm{admissable}\ A,w\notin
A}\mathbb{P}\bigl(\mathcal{C}(0)=A\bigr)\mathbb{P}\bigl(w\stackrel{r}{\leftrightarrow
}y\mid\mathcal{C}(0)=A\bigr),
\end{eqnarray*}
where ``$A$ admissible'' means that $A$ is a connected subgraph
of $T\times T$ containing $0$, $x$, and a path of length $\le r$
between $0$ and $x$. Note that for admissible $A$ with $w\notin A$
we have $\mathbb{P}(w\stackrel{r}{\leftrightarrow}y\mid\mathcal
{C}(0)=A)=\mathbb{P}(w\stackrel{r}{\leftrightarrow}y\mbox{ off }A)$
where the event $\{w\stackrel{r}{\leftrightarrow}y$ off $A\}$ means
that there exists an open path of length at most $r$ connecting $w$
to $y$ which avoids the vertices of $A$. At this point, we can remove
the condition $w\notin A$ since in this case the event $\{w\stackrel
{r}{\leftrightarrow}y$\vadjust{\eject} off~$A\}$ is empty. We get
\[
\mathbb{P}\bigl(0\stackrel{r}{\leftrightarrow}x,  w\stackrel
{r}{\leftrightarrow}y, 0\nleftrightarrow w\bigr)=\sum_{\mathrm{admissable}
\ A}\mathbb{P}\bigl(\mathcal{C}(0)=A\bigr)\mathbb{P}\bigl(w\stackrel{r}{\leftrightarrow
}y\mbox{ off }A\bigr).
\]
Now, obviously
\[
\mathbb{P}\bigl(0\stackrel{r}{\leftrightarrow}x\bigr)\mathbb{P}\bigl(w\stackrel
{r}{\leftrightarrow}y\bigr)=\sum_{\mathrm{admissable}\ A}\mathbb{P}\bigl(\mathcal
{C}(0)=A\bigr)\mathbb{P}\bigl(w\stackrel{r}{\leftrightarrow}y\bigr),
\]
and we subtract these two equalities and get
%
\begin{eqnarray}\label{eq:revBKstp1}
&&\mathbb{P}\bigl(0\stackrel{r}{\leftrightarrow}x,  w\stackrel
{r}{\leftrightarrow}y, 0\nleftrightarrow w\bigr)\nonumber\\
&&\qquad=\mathbb{P}\bigl(0\stackrel
{r}{\leftrightarrow}x\bigr)\mathbb{P}\bigl(w\stackrel{r}{\leftrightarrow}y\bigr)\\
&&\qquad\quad{}-\sum_{\mathrm{admissable}\ A}\mathbb{P}\bigl(\mathcal{C}(0)=A\bigr)\mathbb
{P}\bigl(w\stackrel{r}{\leftrightarrow}y\mbox{ only on }A\bigr),\nonumber
\end{eqnarray}
where the event $\{w\stackrel{r}{\leftrightarrow}y\mbox{ only on }A\}$
means that there exists an open path between $w$ and $y$ of length
at most $r$ and any such path must have a vertex in $A$. Denoting
such a vertex by $v$ we get that $\{w\stackrel{r}{\leftrightarrow}v\}
\circ\{v\stackrel{r}{\leftrightarrow}y\}$.
Hence, for any subgraph $A$ of $T\times T$ we have
\[
\mathbb{P}\bigl(w\stackrel{r}{\leftrightarrow}y\mbox{ only on }A\bigr)\leq\sum
_{v\in A}\mathbb{P}\bigl(\bigl\{w\stackrel{r}{\leftrightarrow}v\bigr\}\circ\bigl\{v\stackrel
{r}{\leftrightarrow}y\bigr\}\bigr).
\]
Putting this into the second term of the right-hand side of (\ref{eq:revBKstp1})
and changing the order of summation gives that we can bound this term
from above by
\begin{eqnarray*}
\mbox{The 2nd  term in
(\ref{eq:revBKstp1})}&\le&
\sum_{v\in T\times T}\mathbb{P}\bigl(\bigl\{v\stackrel{r}{\leftrightarrow
}w\bigr\}\circ\bigl\{v\stackrel{r}{\leftrightarrow}y\bigr\}\bigr)
\!\!\!\sum_{A\ \mathrm{admissible},v\in A}\!\!\!\mathbb{P}\bigl(\mathcal{C}(0)=A\bigr)\\
& = &\sum_{v\in T\times T}\mathbb{P}\bigl(\bigl\{v\stackrel{r}{\leftrightarrow}w\bigr\}
\circ\bigl\{v\stackrel{r}{\leftrightarrow}y\bigr\}\bigr)\mathbb{P}\bigl(0\stackrel
{r}{\leftrightarrow}x, 0\leftrightarrow v\bigr).
\end{eqnarray*}
Now, if $0\stackrel{r}{\leftrightarrow}x$ and $0\leftrightarrow v$
then there exists $u$ such that the events $0\stackrel
{r}{\leftrightarrow}u$,
$u\leftrightarrow v$ and $u\stackrel{r}{\leftrightarrow}x$ occur
disjointly. We use the BK inequality and get
\[
\le\sum_{u,v\in T\times T}\mathbb{P}\bigl(v\stackrel{r}{\leftrightarrow
}w\bigr)\mathbb{P}\bigl(v\stackrel{r}{\leftrightarrow}y\bigr)\mathbb{P}\bigl(0\stackrel
{r}{\leftrightarrow}u\bigr)\mathbb{P}(u\leftrightarrow v)\mathbb
{P}\bigl(u\stackrel{r}{\leftrightarrow}x\bigr).
\]
Let us now sum (\ref{eq:revBKstp1}) over $x$ and $y$ and use the
estimate above for the second term on its right-hand side. We get
%
\begin{eqnarray}\label{eq:offpointa}
&&\sum_{x,y\in T\times T}\mathbb{P}\bigl(0\stackrel{r}{\leftrightarrow}x,
w\stackrel{r}{\leftrightarrow}y, 0\nleftrightarrow
w\bigr)\nonumber\\[-8pt]\\[-8pt]
&&\qquad\ge G(r)^{2}-G(r)^{2}\sum_{u,v}\mathbb{P}\bigl(0\stackrel{r}{\leftrightarrow
}u\bigr)\mathbb{P}(u\leftrightarrow v)\mathbb{P}\bigl(v\stackrel
{r}{\leftrightarrow}w\bigr)\nonumber
\end{eqnarray}
(please record this inequality in this form as we will need it later
in Lemma~\ref{lem:AO1loglog}). With Lem\-ma \ref{lem:openoded}, we
are done.
\end{pf*}
\begin{lem}
\label{lem:AO2}There exists some constant $C_{1}$ such that $G(2r)\ge
r^{-C_{1}}G(r)^{2}$
for all $r\ge2$.
\end{lem}
\begin{pf}
Fix $w$ as in Lemma \ref{lem:AO1_1}---to be more precise, let
$w$ be of minimal distance from $0$ satisfying the conclusion of
Lemma \ref{lem:AO1_1}. The proof uses a~modification argument, namely
we show that by a modification that ``costs'' no more than $r^{C}$,
one gets from the event of Lemma \ref{lem:AO1_1} to the event $\{
0\stackrel{2r+|w|}{\longleftrightarrow}y\}$,
summing over the probabilities of which would allow to lower bound
$G(2r+|w|)$. This is enough because, clearly,
\[
B_{\chem}(2r+|w|)\subset\bigcup_{x\in B_{\chem}(2r)}B(x,|w|)
\]
[note that the ball on the right, $B(x,|w|)$, is in the original
graph and not in the intrinsic metric] and hence $G(2r+|w|)\le
G(2r)\cdot|B(|w|)|\le G(2r)\cdot2d^{|w|}$
and since $|w|\le C\log r$,
\[
G(2r)\ge r^{-C}G(2r+|w|)
\]
so it is enough to lower bound $G(2r+|w|)$. Returning to the modification,
the process would be to take the clusters containing the path from
$0$ to $x$ and from $xw$ to $y$ and connect them by the shortest
possible path. Formally, we do as follows. Let $A$ be the collection
of all triplets $(\pi,x,y)$ such that $x,y\in T\times T$ and $\pi$
is some configuration on $T\times T$ such that in $\pi$ we have
$0\stackrel{r}{\leftrightarrow}x$, $xw\stackrel{r}{\leftrightarrow}y$
and $0\nleftrightarrow y$. Let $B$ be the collection of all couples
$(\pi,y)$~whe\-re~$\pi$ is some configuration such that $0\stackrel
{2r+|w|}{\longleftrightarrow}y$.
We shall construct a $\varphi\dvtx A\to B$ with the following two
properties:
\begin{itemize}
\item$\varphi$ is no more than $r^{C}$ to $1$.
\item The Radon--Nikodym derivative of $\varphi$ is bounded below by $r^{-C}$
(we consider $A$ and $B$ as measure spaces with the counting measure
for $x$ and $y$ and the usual product measure for $\pi$).
\end{itemize}
Clearly, once $\varphi$ is constructed, we would get
\[
G(2r+|w|)=|B|\ge r^{-C}|A|\ge r^{-C}\cdot\tfrac{1}{2}G(r)^{2},
\]
where the last inequality is by Lemma \ref{lem:AO1_1} and where $|A|$
and $|B|$ stand for the total measure of $A$ and $B$, respectively.
So we need only construct $\varphi$.

The construction is as follows. Let $(\pi,x,y)\in A$. Let $\gamma$
be a shortest path from $x$ to $xw$ (choose $\gamma$ arbitrarily,
e.g., first walk on the first tree and then on the second).
Let $e$ be the last point on $\gamma$ which is in $\mathcal{C}(x)$.
Let $f$ be the first point on $\gamma$ after $e$ which is in $\mathcal{C}(xw)$.
Let $\pi'$ be the configuration one gets by opening every edge of
$\pi$ on the piece of $\gamma$ between $e$ and $f$. See Figure
\ref{fig:phi}.%
%
\begin{figure}

\includegraphics{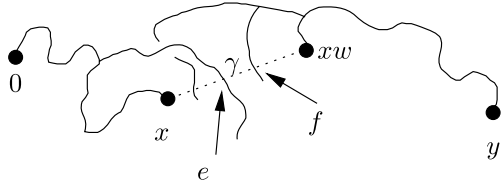}

\caption{The construction of the modification $\varphi$.}
\label{fig:phi}
\end{figure}
Define $\varphi(\pi,x,y)=(\pi',y)$. Clearly, the Radon--Nikodym derivative
is equal to
\[
\biggl(\frac{p}{1-p}\biggr)^{\#\mathrm{closed}\ \mathrm{edges}\ \mathrm{in}\ (e,f)}\ge
\biggl(\frac{p}{1-p}\biggr)^{|\gamma|}\ge\biggl(\frac{p}{1-p}\biggr)^{C\log
r}=r^{-C},
\]
where $p=p_{c}(T\times T)$. Recall that $|\gamma|\le C\log r$ by
Lemma \ref{lem:AO1_1}. To show that $\varphi$ is no more than $r^{C}$
to $1$, examine one couple $(\pi',y)\in B$. If $\varphi(x,y,\pi)=(\pi',y)$,
then all edges between $e$ and $f$ must be pivotal for the connection
$0\leftrightarrow y$. Since there can be no more than $2r+|w|$ edges
which are pivotal for the connection $0\stackrel
{2r+|w|}{\longleftrightarrow}y$,
we see that $e$ has no more than $2r+|w|$ possibilities. Since $|e-x|\le|w|$,
we see that $x$ has no more than $(d-1)^{|w|}\le(d-1)^{C\log r}=r^{C}$
possibilities. Once $x$ is fixed so is $\gamma$. The original configuration
on $\gamma$ has $2^{|w|}\le r^{C}$ possibilities. The shows that
$\varphi$ is no more than $r^{C}$ to 1 and finishes the lemma.
\end{pf}
\begin{lem}
\label{lem:AO3}$G(r)\le Cr^{C}$.
\end{lem}
\begin{pf}
This is a more-or-less direct corollary of Lemmas \ref{lem:Goer}
and \ref{lem:AO2}. Let~$C_{1}$ be the constant from Lemma
\ref{lem:AO2}. Assume by contradiction that for some~$r$, $G(r)>(4r)^{C_{1}}$.
Then by applying Lemma \ref{lem:AO2} repeatedly,
\begin{eqnarray*}
G(2r) & \ge & r^{-C_{1}}G(r)^{2}>(2^{4}r)^{C_{1}},\\
G(4r) & \ge &(2r)^{-C_{1}}G(2r)^{2}>(2^{8-1}r)^{C_{1}},\\
& \vdots &\\
G(2^{k+1}r) & \ge &(2^{k}r)^{-C_{1}}G(2^{k}r)^{2}>(2^{2^{k+2}-\sum
_{l=1}^{k}l2^{k-l}}r)^{C_{1}}>(2^{2^{k+1}}r)^{C_{1}},
\end{eqnarray*}
which means that $G(s)$ increases exponentially in $s$, contradicting
Lemma~\ref{lem:Goer}.
\end{pf}

We now repeat the arguments of Lemmas \ref{lem:AO1_1}--\ref{lem:AO3},
but use Lemma \ref{lem:AO3} as an input to get better results.
\begin{lem}
\label{lem:AO1loglog}There exists some constant $C_{1}$ such that
for any $r\ge3$ and any $w$ with $|w|\ge C_{1}\log\log r$,
\[
\mathbb{E}\bigl(\bigl|\bigl\{(x,y)\dvtx0\stackrel{r}{\leftrightarrow}x,  xw\stackrel
{r}{\leftrightarrow}y, 0\nleftrightarrow
y\bigr\}\bigr|\bigr)\ge\tfrac{1}{2}G(r)^{2}.
\]
\end{lem}
\begin{pf}
We start the calculation from (\ref{eq:offpointa}), which, we recall,
stated that
\begin{eqnarray*}
&&\sum_{x,y\in T\times T}\mathbb{P}\bigl(0\stackrel{r}{\leftrightarrow}x,
w\stackrel{r}{\leftrightarrow}y, 0\nleftrightarrow w\bigr)\\
&&\qquad\geq
G(r)^{2}-G(r)^{2}\sum_{u,v}\mathbb{P}\bigl(0\stackrel{r}{\leftrightarrow
}u\bigr)\mathbb{P}(u\leftrightarrow
v)\mathbb{P}\bigl(v\stackrel{r}{\leftrightarrow}w\bigr).
\end{eqnarray*}
Recall also that the sum on the right-hand side is denoted by $\nabla(w;r)$.
We now separate the sum into two parts, according to whether $\max\{
|u|,|v|\}\le\log^{2}r$
or not. The first case is calculated exactly as in Lemma \ref{lem:AO1_1},
as follows
\begin{eqnarray*}
\hspace*{-4pt}&&\phantom{\mbox{By Lemma \ref{lem:Oded}}}\sum_{|u|,|v|\le\log^{2}r}\mathbb{P}
\bigl(0\stackrel{r}{\leftrightarrow
}u\bigr)\mathbb{P}(u\leftrightarrow v)\mathbb{P}\bigl(v\stackrel
{r}{\leftrightarrow}w\bigr)\\
\hspace*{-4pt}&&
\mbox{By Lemma \ref{lem:Oded}} \qquad\le\sum_{|u|,|v|\le\log
^{2}r}|u|^{2}d(u,v)^{2}d(v,w)^{2}(d-1)^{-(|u|+d(u,v)+d(v,w))/2}\\
\hspace*{-4pt}&&\phantom{\mbox{By Lemma \ref{lem:Oded}}}\qquad\le C|w|^{2}\log^{12}r\sum_{|u|,|v|\le\log
^{2}r}(d-1)^{-(|u|+d(u,v)+d(v,w))/2}\\
\hspace*{-4pt}&&\mbox{\hspace*{29.2pt}By (\ref{eq:yzr})}\qquad \le
C|w|^{2}\log^{12}r\cdot9(\log^{2}r)^{8}(d-1)^{-|w|/2}\\
\hspace*{-4pt}&&\phantom{\mbox{By Lemma \ref{lem:Oded}}}\qquad=C|w|^{2}(d-1)^{-|w|/2}\log^{28}r
\end{eqnarray*}
and this is $\mbox{$\le$}\frac{1}{4}$ if only $C_{1}$ is sufficiently large.
Now assume $\max\{|u|,|v|\}>\log^{2}r$. As $u$ and $v$ are symmetric,
we may assume $|u|>\log^{2}r$. Let $L$ be the level of $T\times T$
which contains $u$, namely $L=\{z\dvtx z_{1}=u_{1},z_{2}=u_{2}\}$. Clearly,
$|L|=(d-1)^{|u|}$. By Lemma \ref{lem:AO3},
\[
\sum_{z\in L}\mathbb{P}\bigl(0\stackrel{r}{\leftrightarrow}z\bigr)\le G(r)\le Cr^{C}.
\]
But $L$ is completely symmetric, so for any $z\in L$, $\mathbb
{P}(0\stackrel{r}{\leftrightarrow}z)=\mathbb{P}(0\stackrel
{r}{\leftrightarrow}u)$.
We get
\[
\mathbb{P}\bigl(0\stackrel{r}{\leftrightarrow}u\bigr)=\frac{1}{|L|}\sum_{z\in
L}\mathbb{P}\bigl(0\stackrel{r}{\leftrightarrow}z\bigr)\le\frac{Cr^{C}}{(d-1)^{|u|}}.
\]
We need to compare this estimate to the estimate of (\ref{eq:yzr})
from Lemma \ref{lem:AO1_1} which we also used above, namely to $(d-1)^{-|u|/2}$.
So we write this as
\[
\mathbb{P}\bigl(0\stackrel{r}{\leftrightarrow}u\bigr)\le(d-1)^{-|u|/2}\cdot\frac
{Cr^{C}}{(d-1)^{(\log^{2}r)/2}}\le(d-1)^{-|u|/2}\cdot Cr^{-15}.
\]
This allows us to write
\begin{eqnarray*}
&&\mathop{\sum_{u,v}}_{|u|>\log^{2}r}
\mathbb{P}\bigl(0\stackrel{r}{\leftrightarrow}u\bigr)\mathbb{P}(u\leftrightarrow
v)\mathbb{P}\bigl(v\stackrel{r}{\leftrightarrow}w\bigr)\\
&&\qquad\le Cr^{-15}\cdot Cr^{4}\sum_{|u|,|v|\le
r}(d-1)^{-(|u|+d(u,v)+d(v,w))/2}\\
&&\qquad\stackrel{\fontsize{8.36pt}{8.36pt}\selectfont{(\ref{eq:yzr})}}{\le
}Cr^{-11}\cdot Cr^{8}(d-1)^{-|w|/2}
\end{eqnarray*}
and we see that this part of the sum is in fact negligible (if $r$
is sufficiently large or if $C_{1}$ is chosen sufficiently large).
This shows that $\nabla(w;r)\le\frac{1}{2}$ and concludes the lemma.
\end{pf}
\begin{lem}
\label{lem:AO2loglog}$G(r)\le Cr(\log r)^{C}$.
\end{lem}
\begin{pf}
This is nothing more than repeating the arguments of Lemmas~\ref{lem:AO2}
and~\ref{lem:AO3}. Let us verify some of the details. We first show
%
\begin{equation}\label{eq:AO2loglog}
G(2r)\ge\frac{G(r)^{2}}{r\log^{C}r},
\end{equation}
that is, the analog of Lemma \ref{lem:AO2}. We again construct a
$\varphi\dvtx A\to B$
($A$ and $B$ being exactly as in Lemma \ref{lem:AO2}) which is
no more than $r\log^{C}r$ to 1, and with the Radon--Nikodym derivative
bounded below by $\log^{-C}r$. The construction of $\varphi$ is
identical, that is, we take the shortest path $\gamma$ from $x$ to
$xw$, let $e$ be the last point of $\mathcal{C}(x)$ on $\gamma$
and $f$ be the first point of $\mathcal{C}(xw)$ on $\gamma$ after
$e$. This time, though, because $|w|\approx\log\log r$, the Radon--Nikodym
derivative of~$\varphi$ will be $\ge(p/(1-p))^{|w|}\ge\log^{-C}r$.
To invert~$\varphi$, we need to find~$e$, $x$ and~$\omega$. $e$
still has $2r+|w|$ possibilities, but given $e$, $x$ has only
$(d-1)^{|w|}\le\log^{C}r$
possibilities, and $\omega$ has only $2^{|w|}\le\log^{C}r$ possibilities.
This shows (\ref{eq:AO2loglog}).

Concluding from (\ref{eq:AO2loglog}) the lemma is identical to the
proof of Lemma \ref{lem:AO3}---if $G(r)>r(4\log r)^{C}$ for some
$r$ then it starts growing exponentially---and we will omit it.
\end{pf}
\begin{lem}
\label{lem:ext}In the extrinsic metric,
\[
\mathbb{E}\bigl(\mathcal{C}(0)\cap B(r)\bigr)\le Cr^{3}.
\]
\end{lem}

The idea is quite simple. We consider $\mathcal{C}(0)$ as a branching
process embedded into $T\times T$ and show that it escapes from $0$
with ``positive speed.'' Our ``time'' for the branching process
is the intrinsic distance $d_{\chem}(0,x)$ so escaping in positive
speed mean simply that $d_{\chem}(0,x)\approx|x|$ (this is step 5
of the sketch on page \pageref{sub:Step-5}). Here, are the details.
It is enough to prove:
\begin{lem}
\label{lem:exiclm}$\mathbb{E}((B_{\chem}(4r^{2})\setminus
B_{\chem}(r^{2}))\cap B(r))\le Cr^{-2}$.\vadjust{\eject}
\end{lem}
\begin{pf*}{Proof of Lemma \ref{lem:ext} given Lemma \ref{lem:exiclm}}
We apply
Lemma \ref{lem:exiclm} with the parameter $r_{\mathrm{Lemma}\ \mbox{\fontsize{8.36pt}{8.36pt}\selectfont{\ref
{lem:exiclm}}}}=r2^{k}$
and get
\begin{eqnarray*}
&&\mathbb{E}\bigl(\bigl(B_{\chem}(4^{k+1}r^{2})\setminus
B_{\chem}(4^{k}r^{2})\bigr)\cap
B(r)\bigr)\\
&&\qquad\le\mathbb{E}\bigl(\bigl(B_{\chem}(4^{k+1}r^{2})\setminus
B_{\chem}(4^{k}r^{2})\bigr)\cap
B(r2^{k})\bigr)\le\frac{C}{r^{2}4^{k}},
\end{eqnarray*}
which we sum over $k$ to get
\[
\mathbb{E}\bigl(\bigl(\mathcal{C}(0)\setminus B_{\chem}(r^{2})
\bigr)\cap B(r)\bigr)\le\frac{C}{r^{2}}.
\]
For the interior part, we just use Lemma \ref{lem:AO2loglog} and get
\[
\mathbb{E}\bigl(B_{\chem}(r^{2})\cap B(r)\bigr)\le\mathbb{E}
(B_{\chem}(r^{2}))\le Cr^{2}\log^{C}r
\]
proving the claim.
\end{pf*}
\begin{pf*}{Proof of Lemma \ref{lem:exiclm}} Clearly, we may assume $r$
is sufficiently
large. We shall define a sequence of subsets of $\mathcal{C}(0)$,
$\{ \partial_{m}\} _{m=1}^{\infty}$ (here $\partial_{m}$
is just a~notation). Intuitively, you should consider $\partial_{m}$
simply as $\partial B_{\chem}(mr)$---I~could not make the proof
work with this definition of $\partial_{m}$ so a somewhat more complicated,
inductive definition will be used [formally also $\partial_{m}$ is
a~set of vertices while $\partial B_{\chem}(mr)$ is a set of edges].
It will preserve some of the feeling of $\partial B_{\chem}(mr)$
as every $v\in\partial_{m}$ will satisfy $\frac{1}{2}mr\le d_{\chem
}(0,v)\le mr$,
and since for any $x\in\mathcal{C}(0)$ there will be a path from
$0$ to $x$ visiting each $\partial_{m}$ in turn and spending $\mbox{$\approx$}r$
steps between any two levels. This path will not be a~geodesic in
$\mathcal{C}(0)$ (i.e., a shortest possible path), but this is not
important.

During the induction process, we shall expose parts of $\mathcal{C}(0)$.
We shall denote the exposed part by $E_{m}$ [which is formally a
collection of edges of~$\mathcal{C}(0)$ and of $\partial\mathcal{C}(0)$,
though we shall say about a vertex $v$ that it is ``in~$E_{m}$''
if there exists some path of open edges in $E_{m}$ from $0$ to $v$].
$\partial_{m}$ will be the ``boundary of $E_{m}$'' in the sense
that any vertex of $\partial_{m}$ is in $E_{m}$, and any vertex
in $E_{m}\setminus\partial_{m}$ is fully exposed, that is, all edges
coming out of it (whether open or closed) are in $E_{m}$. We start
with $\partial_{0}=\{0\}$ and $E_{0}$ having no edges. Note that
saying about $E_{m}$ that it is ``exposed'' is not just a name,
it carries some meaning, namely that for any $m$ and any set of edges
$A$, one can infer whether $E_{m}=A$ or not merely by examining
the states of the edges of $A$. We will keep this property through
the induction.\vspace*{10pt}

1. \textit{The construction of $\partial_{m}$.} Assume $E_{m}$ has
already been calculated. Let $y\in\partial_{m}$. We will now construct
the ``children of $y$'' which will belong to $\partial_{m+1}$.
For this purpose, let $Q=Q(y)$ be the set of all vertices $q$ satisfying
that $q$ is connected to $y$ by an open path of length $\le r$
off $E_{m}$ (we are making the exception that $y$ is in $E_{m}$,
but no other vertex of the path can be in~$E_{m}$, this is the precise
meaning of ``off'' here). Clearly, $\mathbb{E}(Q)\le G(r)$ and
by Lemma~\ref{eq:AO2loglog} we get
%
\begin{equation}\label{eq:Qrlogr}
\mathbb{E}(Q)\le C_{1}r(\log r)^{C}.
\end{equation}
We make at this point the convention that every formula involving
$Q$ is in fact conditioned over $E_{m}$. For example, (\ref{eq:Qrlogr})
should be read as $\mathbb{E}(Q |  E_{m})\le C_{1}r(\log r)^{C}$.
We need to examine two special parts of $Q$---the vertices ``close
to $y$'' and the vertices ``beyond $y$.'' For the first part,
let $\ell$ be defined by
\[
|B(\ell)|\le\sqrt{r}
\]
and $\ell$ being the maximal with this property (note that this is
the usual ball in our original graph $T\times T$, and we have $\ell
\approx\log r$).
Clearly,
%
\begin{equation}\label{eq:Qclose}
\mathbb{E}\bigl(|Q\cap B(y,\ell)|\bigr)\le|B(y,\ell)|\le\sqrt{r}.
\end{equation}
For the second part, fix some $k_{1}$ and $k_{2}$, and let $L$
be the corresponding level, that is,
\[
L=\{z\dvtx d(z_{i},y_{i})=k_{i},  i=1,2\}.
\]
A straightforward calculation shows that
\[
\bigl|\{z\in L\dvtx |z|<|y|\}\bigr|\le|L|e^{-c|k|}.
\]
However, $L$ is completely symmetric with respect to $y$. Therefore,
\[
\mathbb{E}\bigl|\bigl\{ z\in L\dvtx z\stackrel{r}{\leftrightarrow}y\mbox{ and
}|z|<|y|\bigr\}\bigr|\le G(r)\frac{|\{ z\in L\dvtx |z|<|y|\}
|}{|L|}\le G(r)e^{-c|k|}.
\]
Requiring that the connection is off $E_{m}$ only makes things worse,
so
%
\begin{equation}\label{eq:before gen k}
\mathbb{E}\bigl|\{ z\in Q\cap L\dvtx |z|<|y|\}\bigr|\le
G(r)e^{-c|k|}.
\end{equation}
Summing this over all $|k|>\ell$ gives
%
\begin{equation}\label{eq:Qbefore}
\mathbb{E}\bigl|\{z\in Q\setminus B(y,\ell)\dvtx|z|<|y|\}\bigr|\le\sum_{|k|>\ell
}G(r)e^{-c|k|}\le Ce^{-c\ell}G(r)\le Cr^{1-c}.\hspace*{-28pt}
\end{equation}
We may sum up (\ref{eq:Qclose}) and (\ref{eq:Qbefore}) as
%
\begin{equation}\label{eq:Qcb}
\mathbb{E}\bigl|\{z\in Q\dvtx|z|<|y|+\ell\}\bigr|\le C_{2}r^{1-c}.
\end{equation}

Equipped with the estimates (\ref{eq:Qrlogr}) and (\ref{eq:Qcb})
we may now proceed to define an important element of the construction,
the parameter $s$. Define $Q_{s}=Q_{s}(y)$ similarly to $Q$ but
with the requirement that the shortest open path (off $E_{m}$) has
length exactly $s$ [below we will denote this by ``$d_{\mathrm{off}\ E_{m}}(y,z)=s$''
for short]. We get $Q=\dotcup_{s=1}^{r}Q_{s}$. This means that for
at least $\frac{2}{3}$ of the $s$ between $\frac{1}{2}r$ and $r$
we must have
%
\begin{equation}\label{eq:sgood0}
\mathbb{E}|Q_{s}|\le6C_{1}(\log r)^{C}
\end{equation}
and similarly,
%
\begin{equation}\label{eq:sgood1}
\mathbb{E}\bigl|\{z\in Q_{s}\dvtx|z|<|y|+\ell\}\bigr|\le6C_{2}r^{-c}.
\end{equation}
Of course, the set of $s$ which satisfy (\ref{eq:sgood0}) may be
different than the set of~$s$ which satisfy (\ref{eq:sgood1}), but
they intersect and we choose one $s$ which satisfies both [sometimes
we will denote it by $s(y)$ for clarity]. We note that the set of
such ``good'' $s$ is a random set which depends on $E_{m}$ and
on $y$, and this is exactly why we cannot define $\partial_{m}=\partial
B_{\chem}(0,mr)$.

We may now complete the description of the construction. We
define~$E_{m+1}$ to be the set of all edges $(z,x)$ for all $z$ satisfying
that there exists a~$y\in\partial_{m}$ such that $d_{\mathrm{off}\
E_{m}}(y,z)<s(y)$;
and all $x$ which are a neighbor of $z$. We define $\partial_{m+1}$
to be the boundary of $E_{m+1}$ in the sense above. We get that for
any $z\in\partial_{m+1}$ there exists some $y\in\partial_{m}$ such
that the $d_{\mathrm{off}\ E_{m}}(y,z)=s(y)$ (but not vice versa---see
Figure \ref{fig:swallow}).

%
\begin{figure}

\includegraphics{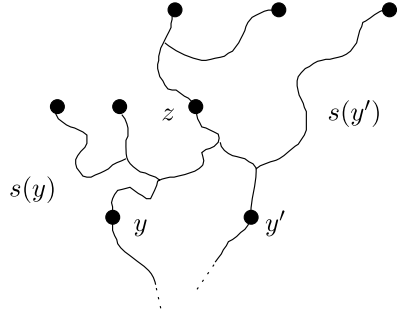}
\vspace*{-4pt}

\caption{The point $z$ has distance exactly $s(y)$ from
$y$ but is not included in $\partial_{m+1}$ because its distance
from $y'$ is less than $s(y')$.}
\label{fig:swallow}
\vspace*{-4pt}
\end{figure}

It will be convenient to note at this point that
%
\begin{equation}\label{eq:BchemE}
B_{\chem}\bigl(\tfrac{1}{2}mr\bigr)\subset E_{m}\subset B_{\chem
}(mr)
\end{equation}
(here $E_{m}$ is considered as a set of vertices). This is a simple
consequence of the fact that all $s(y)$ for all $y$ are in $[\frac{1}{2}r,r]$.
Therefore the direction $E_{m}\subset B_{\chem}(mr)$ is immediate. For
the other direction let $x$ satisfy that $d_{\chem}(x,E_{m})\le\frac{1}{2}r$.
Let $y\in\partial_{m}$ be the closest point in $E_{m}$ to $x$.
Then there exists an open path of length $\le\frac{1}{2}r<s(y)$ off
$E_{m}$ from $y$ to $x$. Hence, $x\in E_{m+1}$ by definition. Hence,
\begin{eqnarray*}
B_{\chem}\bigl(\tfrac{1}{2}(m+1)r\bigr) & = & \bigl\{ x\dvtx d_{\chem}
\bigl(x,B_{\chem}\bigl(\tfrac{1}{2}mr\bigr)\bigr)\le\tfrac{1}{2}r\bigr\} \\
\mbox{inductively\qquad} & \subset &\bigl\{x\dvtx d_{\chem}(x,E_{m})\le\tfrac{1}{2}r\bigr\}\\
\mbox{by the argument above\qquad} &\subset& E_{m+1}
\end{eqnarray*}
showing (\ref{eq:BchemE}).

2. \textit{The number of bad paths.} Let $x\in\partial_{m}$ be some
vertex. Then there exists a~sequence $0=x_{0},\ldots,x_{m}=x$ with
$x_{i}\in\partial_{i}$ and such that for each\vadjust{\eject} $i=0,\ldots,m-1$ we
have that $d_{\mathrm{off}\ E_{i}}(x_{i},x_{i+1})=s(x_{i})$. This
sequence might not be unique but this is not important. We call such
sequences ``$\partial$-paths.''
\begin{defn}
We say that a given $\partial$-path $x_{0},\ldots,x_{m}$ is bad
if there are $>m/\sqrt{\log r}$ values of $i\in\{0,\ldots,m-1\}$
such that $|x_{i+1}|<|x_{i}|+\ell$.
\end{defn}

Let us estimate the expected number of bad $\partial$-paths. Let
$I\subset\{0,\ldots,m-1\}$ be some set of indices with $|I|>m/\sqrt
{\log r}$.
Let $E(I;m)$ be the expected number of $\partial$-paths such that
$|x_{i+1}|<|x_{i}|+\ell$ for every $i\in I$. We get, directly from
the choice of $s$ above,
%
\begin{equation}\label{eq:EImrecurs}
E(I;m)\le\cases{
E(I;m-1)\cdot C(\log r)^{C}, &\quad $m-1\notin I$,\cr
E(I\setminus\{m-1\};m-1)\cdot Cr^{-c}, &\quad $m-1\in I$,}
\end{equation}
[seeing (\ref{eq:EImrecurs}) is a standard exercise in the ``off
method''---one conditions on~$E_{m-1}$ and examines each path
individually. Because being connected off~$E_{m-1}$ does not examine
the edges of~$E_{m-1}$ at all, the estimates (\ref{eq:sgood0}) and~(\ref{eq:sgood1})
hold after the conditioning. Here is where we need
the property that $E_{m}$ is ``exposed,'' i.e., that conditioning
over $E_{m}=A$ gives no information on edges not in $A$]. We apply~(\ref{eq:EImrecurs}) recursively and get
\[
E(I;m) \le(C(\log r)^{C})^{m-|I|}(Cr^{-c}
)^{|I|}
\]
and since $|I|>m/\sqrt{\log r}$,
\[
\le\exp\bigl(m
\bigl(C\log\log r-c\sqrt{\log r}\bigr)\bigr).
\]
Summing over all $I$ can be bounded roughly by multiplying this by
$2^{m}$, and we end with
%
\begin{equation}\label{eq:part2}
\mathbb{E}|\{\mbox{bad }\partial\mbox{-paths}\}|\le\exp\bigl(m
\bigl(C-c\sqrt{\log r}\bigr)\bigr)
\end{equation}
and in particular we get the same estimate for the probability that
even one bad $\partial$-path exists.\vspace*{10pt}

3. \textit{The contribution of the bad part.} Examine a good path
$x_{0},\ldots,x_{m}$. We know it contains $<m/\sqrt{\log r}$ ``bad''
$i$ for which $|x_{i+1}|<|x_{i}|+\ell$. We still need to estimate
how much ``damage'' can each such bad $i$ cause, that is, upper
bound $|x_{i}|-|x_{i+1}|$. For this purpose, we return to (\ref
{eq:before gen k}),
which shows that
\[
\mathbb{E}\bigl|\{z\in Q\dvtx|z|<|y|-j\}\bigr|\le\sum_{|k|>j}G(r)e^{-c|k|}\le
G(r)e^{-cj} \qquad\forall y\in\partial_{i}\ \forall i
\]
and in particular there exists some constant $C_{3}$ such that
%
\begin{equation}\label{eq:defC3}
\mathbb{E}\bigl|\{z\in Q\dvtx|z|<|y|-C_{3}\log r\}\bigr|\le r^{-10} \qquad\forall y\in
\partial_{i}\ \forall i.
\end{equation}
We now wish to sum over all $v$, so denote by $\mathbf{BE}$ (for
``bad edges'') the set of all $(y,z)$ with $y\in\partial_{i}$,
$i=0,\ldots,m-1$, and $z\in Q_{s(y)}(y)$ such that $|z|<|y|-C_{3}\log r$.
So we want to estimate $\mathbb{P}(\mathbf{BE}\ne\varnothing)$. Recall\vadjust{\eject}
that $\partial_{i}\subset B_{\chem}(ir)$~(\ref{eq:BchemE}). Recall
also our convention about $Q$ which stated that (\ref{eq:defC3})
in fact holds also after conditioning over $E_{i}$. We write
\begin{eqnarray*}
\mathbb{P}(\mathbf{BE}\ne\varnothing) &\le& \mathbb{E}(|\mathbf{BE}|)\le
\sum_{i=0}^{m-1}\sum_{y,z}\mathbb{E}\mathbb{P}\bigl((y,z)\in\mathbf{BE} |
E_{i}\bigr)\\
&\le& \sum_{i=0}^{m-1}\sum_{y}\mathbb{E}\bigl(r^{-10}\mathbb{P}(y\in\partial
_{i} |  E_{i})\bigr)\\
&=& r^{-10}\sum_{i=0}^{m-1}\mathbb{E}|\partial_{i}|\le r^{-10}G(mr)\\
&\le&
mr^{-9}(\log mr)^{C}.
\end{eqnarray*}

\vspace*{10pt}

4. \textit{Wrapping it all up.} Let $m\in\{r,\ldots,8r\}$. Our analysis
of bad $\partial$-paths concluded with (\ref{eq:part2}) which states
that
\begin{eqnarray*}
&&\mathbb{P}(\exists\mbox{ bad }\partial\mbox{-path for some
}m\in\{r,\ldots,8r\})\\
&&\qquad\le(7r+1)\exp\bigl(r\bigl(C-c\sqrt{\log r}\bigr)\bigr)\le Cr^{-7}.
\end{eqnarray*}
Denote this event by $\mathcal{B}_{1}$. From part 3 of the proof,
we have that, with probability $\mbox{$\le$} Cr^{-7}$, there exists some
$\partial$-path
and some $i$ such that $|x_{i+1}|\le|x_{i}|-C_{3}\log r$. Denote
this event by $\mathcal{B}_{2}$ and $\mathcal{B}=\mathcal{B}_{1}\cup
\mathcal{B}_{2}$.
If $\neg\mathcal{B}$ happened (the ``good'' case), then for every
$x\in\partial_{m}$ we have
\[
|x|\ge\ell m-\frac{m}{\sqrt{\log r}}\cdot C_{3}\log r\ge cr\log r
\]
for $r$ sufficiently large.

We now move from $x\in\partial_{m}$ to some arbitrary $x\in B_{\chem
}(4r^{2})\setminus B_{\chem}(r^{2})$.
Let~$m_{0}$ be the first $m$ such that $x\in E_{m+1}$, and let
$y\in\partial_{m}$ be the closest point to $x$, so that $d_{\chem
}(y,x)\le r$.
Recall (\ref{eq:BchemE}) which stated that $B_{\chem}(\frac
{1}{2}mr)\subset E_{m}\subset B_{\chem}(mr)$.
We get that $m_{0}\in\{r,\ldots,8r\}$, and therefore (still assuming~$\neg\mathcal{B}$)
$|y|\ge cr\log r$ so
\[
|x|>|y|-d(y,x)\ge cr\log r-r.
\]
This finishes the lemma, as we see that $(B_{\chem
}(4r^{2})\setminus B_{\chem}(r^{2}))\cap B(r)$
can be nonempty only if the bad event $\mathcal{B}$ happened. On
the one hand, by Lemma~\ref{lem:diagram} we see that
\[
\mathbb{E}\bigl(|B_{\chem}(4r^{2})|\cdot\mathbf{1}\{ |B_{\chem
}(4r^{2})|>C(\log r)G(4r^{2})^{2}\} \bigr)\le Cr^{-2}
\]
and by Lemma \ref{lem:AO2loglog} $C(\log r)G(4r^{2})^{2}\le
C(\log r)^{C}r^{4}\le Cr^{5}$.
We get
\begin{eqnarray*}
&&\mathbb{E}\bigl(\bigl(B_{\chem}(4r^{2})\setminus B_{\chem}(r^{2})
\bigr)\cap B(r)\bigr)\\
&&\qquad  \le\mathbb{E}\bigl(|B_{\chem}(4r^{2})\cdot\mathbf{1}\{|B_{\chem
}(4r^{2})|>Cr^{5}\}\bigr)+Cr^{5}\cdot\mathbb{P}(\mathcal{B})\\
&&\qquad\le
Cr^{-2},
\end{eqnarray*}
which proves Lemma \ref{lem:exiclm} and hence Lemma \ref{lem:ext}.
\end{pf*}
\begin{pf*}{Proof of Theorem \ref{theo1}}
By Lemma \ref{lem:ext} and the fact that every
$u$ has $(d-1)^{|u|}$ clones,
%
\begin{equation}\label{eq:2|u|}
\mathbb{P}(0\leftrightarrow u)\le C|u|^{3}(d-1)^{|u|}.
\end{equation}
With this we write
\begin{eqnarray*}
\nabla &=& \sum_{u,v}\mathbb{P}(0\leftrightarrow u)\mathbb
{P}(u\leftrightarrow v)\mathbb{P}(v\leftrightarrow0)\\
\mbox{by (\ref{eq:2|u|})\qquad} &\le& \sum
_{u,v}|u|^{3}d(u,v)^{3}|v|^{3}(d-1)^{-|u|-d(u,v)-|v|}\\
\mbox{as in Lemma \ref{lem:openoded}\qquad} &\le&\biggl(C\sum
_{u_{1},v_{1}}|u_{1}|^{3}d(u_{1},v_{1})^{3}|v_{1}|^{3}(d-1)^{-|u_{1}|-d(u_{1},v_{1})-|v_{1}|}
\biggr)^{2},
\end{eqnarray*}
where the last sum is over $u_{1}$ and $v_{1}$ in the tree $T$
(and not in $T\times T$). Denote by $a$ the point where the paths
from $0$ to $u_{1}$ and from $0$ to $v_{1}$ split. Then
$|u_{1}|=|a|+d(a,u_{1})$,
$d(u_{1},v_{1})=d(u_{1},a)+d(a,v_{1})$ and $|v_{1}|=|a|+d(a,v_{1})$
so all in all
\begin{eqnarray*}
&&\sum_{u_{1},v_{1}}(|u_{1}|d(u_{1},v_{1})|v_{1}|
)^{3}(d-1)^{-|u_{1}|-d(u_{1},v_{1})-|v_{1}|}\\
&&\qquad \le\sum_{u_{1},v_{1}}(|u_{1}|d(u_{1},v_{1})|v_{1}|
)^{3}(d-1)^{-2|a|-2d(a,v_{1})-2d(a,v_{1})}\\
&&\qquad \le C\sum_{u_{1},v_{1}}(|a|d(a,u_{1})d(a,v_{1})
)^{6}(d-1)^{-2|a|-2d(a,v_{1})-2d(a,v_{1})}.
\end{eqnarray*}
Fixing $a$ we may sum over $v_{1}$ in the subtree of $a$ and get
$\sum d(a,v_{1})^{6}(d-1)^{-2d(a,v_{1})}$ which is finite and independent
of $a$. The sum over the $u_{1}$ in the subtree of $a$ gives an
identical contribution. Finally, we sum over $a$ and get a third contribution
identical to the previous two. Hence, this sum is finite and so is
$\nabla$.\vspace*{6pt}
\end{pf*}

\section*{\texorpdfstring{Acknowledgments.}{Acknowledgments}}

I was fortunate enough to hear the proof of Schramm's lemma in a talk
given by Yuval Peres in an AIM workshop on ``Percolation on transitive
graphs,'' May 2008. I wish to thank Yuval Peres for the inspiring
talk and, of course, the hospitality of the American Institute of
Mathematics. Also, I wish to thank Remco van der Hofstad for answering
all my questions on lace expansion, Markus Heydenreich for pointing
out some errors in a~draft version and Johan Tykesson for pointing
me to \cite{AV08}.

%

%
\printaddresses

\end{document}